\documentclass[12pt,a4paper, twoside,reqno]{amsart}
\usepackage{tikz}
\usepackage{amsmath}
\usepackage{amssymb,amsfonts,mathrsfs}
\usepackage{amscd,amssymb,verbatim}
\usepackage[colorlinks=true,linkcolor=blue,citecolor=blue]{hyperref}
\usepackage{mllocal}

\theoremstyle{plain}

\numberwithin{equation}{section}

\begin{document}

\title[Refined Analytic Torsion on Manifolds with Boundary]
{Refined analytic torsion as analytic function on the representation variety and applications}

\author{Maxim Braverman}
\address{Department of Mathematics,
Northeastern University,}\address{
Boston, MA 02115,
USA}

\email{maximbraverman@neu.edu}
\urladdr{www.math.neu.edu/~braverman/}

\author{Boris Vertman}
\address{Mathematisches Institut,
Universit\"at Bonn,}\address{
53115 Bonn,
Germany}
\email{vertman@math.uni-bonn.de}
\urladdr{www.math.uni-bonn.de/people/vertman}

\date{This document was compiled on: \today}

\begin{abstract}
{We prove that refined analytic torsion on a manifold with boundary
is a weakly holomorphic section of the determinant line bundle over the representation variety.
As a fundamental application we establish a gluing formula for refined analytic torsion 
on connected components of the complex representation 
space which contain a unitary point. Finally we provide a new proof of Br\"uning-Ma gluing formula for the Ray-Singer torsion associated to a non-Hermitian connection. Our proof is quite different from the one given by Br\"uning and Ma and uses a temporal gauge transformation.}
\end{abstract}

\maketitle

\tableofcontents

\section{Introduction and statement of the main results}\label{intro}

The Ray-Singer conjecture has been formulated in the seminal paper of Ray and Singer \cite{RS}
and proved independently by Cheeger \cite{Che} and M\"uller \cite{Mue1} for unitary representations. Its importance
stems from the fact that as in the Atiyah-Singer index theorem, it equates analytic with combinatorial 
quantities, the analytic Ray-Singer and the combinatorial Reidemeister torsions. 

By construction, both the analytic Ray-Singer and the combinatorial Reidemeister torsions provide canonical norms on the determinant line of cohomology.
There have been various approaches to obtain a canonical construction of analytic and Reidemeister torsions as elements instead of norms
of the determinant line of the cohomology. These constructions seek to refine the notion of analytic and Reidemeister torsion 
norms on that determinant line, which basically corresponds to fixing a complex phase in the family of complex vectors of length one.

In case of the Reidemeister torsion this has been done by Farber and Turaev \cite{FarTur:absolute} and \cite{FarTur:poincare}.
Refinement of analytic torsion has been studied by the first author jointly with Kappeler in \cite{BraKap:refined} and \cite{BraKap:refined-det}, 
as well as by Burghelea and Haller in \cite{BurHal:complex1} and \cite{BurHal:complex2}. Both notions have subsequently 
been compared by the first author jointly with Kappeler in \cite{BraKap:comparison}. An extension of refined analytic torsion to manifolds with boundary
has been undertaken by the second author \cite{Vertman09} and Lee and Huang \cite{LeeHua:refined} in two different independent constructions. 
Recently, Lee and Huang also compared the two notions of refined analytic torsion on manifolds with boundary in \cite{LeeHua:comparison}.

The fundamental property of the refined analytic and Farber-Turaev torsions is that they define weakly holomorphic  
functions on the complex representation space (see \cite[p.~148]{Gunning70-AnVarieties1} for definition of a weakly holomorphic function). In particular, its restriction to the regular part of the representation space is holomorphic. The main purpose  of the present discussion is an extension of this result
to the refined analytic torsion on manifolds with boundary, introduced by the second author in \cite{Vertman09}. 
As a consequence we establish the gluing property of refined analytic torsion on connected components of  the representation variety that contain a unitary point. 

The gluing formula for refined analytic torsion may be used to prove a gluing result for the Ray-Singer torsion norm for certain non-unitary representations path-connected to a unitary element. However we chose to 
devote the final two sections of the present paper to an alternative proof of the gluing property for the Ray Singer analytic 
torsion for non-unitary representations, which is stronger since we do not single out connected components without unitary elements.

The Ray-Singer theorem has been extended to unimodular representations by M\"uller \cite{Mue2}.
In case of a general non-unitary representation, the quotient of the analytic and Reidemeister torsion norms
admits additional correction terms which have been studied by Bismut and Zhang in \cite{BisZha}. In a separate discussion 
\cite{BraVer}, the authors employ analyticity of refined analytic torsion to provide 
an alternative derivation of the Bismut-Zhang correction terms for the connected components of 
unitary points in the representation variety. 

Both the analytic and combinatorial Reidemeister torsions make sense on compact manifolds 
with boundary after posing relative or absolute boundary conditions. The gluing property of 
the analytic torsion, which is foremost a spectral invariant, is striking and has been proved 
by L\"uck \cite{Lue}, Vishik \cite{Vis} and generalized by Lesch \cite{Les}, under the assumption 
of product metric structures and unitary representations. 

The anomaly of analytic torsion on a compact manifold with boundary with a non-unitary representation 
and general metric structures near the boundary, has been studied by Br\"uning and Ma in \cite{BruMa}.
Recently, Br\"uning and Ma established in a follow-up paper \cite{BruMa:gluing} the gluing formula for analytic torsion
in case of a non-unitary representation, generalizing previous results in \cite{Lue} and \cite{Vis}. 
We present an alternative proof of a result by Br\"uning-Ma \cite{BruMa:gluing}
by a temporal gauge transformation argument.

\section{Refined analytic torsion on manifolds with boundary}\label{explicit-unitary} \

This section reviews the construction by the second author \cite{Vertman09}. \\[-7mm]

\subsection{The flat vector bundle induced by a representation.}
Let $(M^m, g)$ be a compact oriented odd-dimensional
Riemannian manifold with boundary $\partial M$.
Consider a complex representation $\alpha$ of the fundamental group $\pi_1=\pi_1(M)$ on $\C^n$. 
Let $(E_\alpha, \D_\alpha, h^E_\alpha)$ be the induced flat complex vector bundle over $M$ with 
monodromy equal to $\alpha$ and no canonical choice of $h^E_\alpha$ in case $\alpha$ is not unitary.

The flat covariant derivative $\D_\alpha$ acts on sections $\Gamma (E_\alpha)$ and extends by Leibniz rule 
to a twisted differential on $E_\alpha$-valued differential forms $\Omega^*_0(M,E_\alpha)$, where 
the lower index refers to compact support in the open interior of $M$. This defines the 
twisted de Rham complex $(\Omega^*_0(M,E_\alpha), \D_\alpha)$. The metrics $(g, h^E_\alpha)$ induce 
an $L^2$-inner product on $\Omega^*_0(M,E_\alpha)$. We denote the $L^2-$completion of 
$\Omega^*_0(M,E_\alpha)$ by $L^2_*(M,E_\alpha)$. 

Throughout this section the representation $\alpha$ is fixed and we omit the lower index $\alpha$
in the notation of $(E_\alpha, \D_\alpha, h^E_\alpha)$ in most of the discussion. 

Next we introduce the notion of the dual covariant derivative $\D'$. It is defined by requiring
for all $u,v \in \Gamma(E)$ and $X \in \Gamma(TM)$
\begin{align}\label{dual-connection}
dh^E(u,v)[X]=h^E(\D_Xu,v)+h^E(u,\D'_Xv).
\end{align}
In the special case that $\alp$ is unitary, the dual $\D'$ and the original covariant derivative $\D$ coincide. 
As before, the dual $\D'$ gives rise to a twisted de Rham complex $(\Omega^*_0(M,E), \D')$. 

\subsection{Hilbert complexes}
For any differential operator $P$ acting on $\Omega^*_0(M,E)$, we denote by $P_{\min}$
its minimal graph-closed extension in $L^2_*(M,E)$. The maximal closed extension is defined by $P_{\max} := (P^t_{\min})^*$.
By Br\"uning and Lesch \cite[Lemma 3.1]{BL1}, the extensions define Hilbert complexes
$(\domr, \Dr)$, where $\domr:=\dom (\Dr)$, and $(\doma, \Da)$, where $\doma:=\dom (\Da)$. 
The Laplace operators, associated to these Hilbert complexes are respectively defined as
\begin{align*}
\triangle_{\textup{rel}}:&=\Dr^*\Dr+\Dr\Dr^*, \\
\triangle_{\textup{abs}}:&=\Da^*\Da+\Da\Da^*.
\end{align*}
Similar definitions hold for the dual connection $\nabla'$ and for 
the Laplace operators $\triangle'_{\textup{rel}}$ and $\triangle'_{\textup{abs}}$ 
of the Hilbert complexes $(\domr', \Dr')$ and $(\doma', \Da')$ respectively.  
The difference $(\D-\D')$ is a bounded endomorphism valued operator and hence the equality of domains
\begin{align}\label{equality}
\domr=\domr', \quad \doma=\doma'.
\end{align}

The following theorem, compare \cite[Theorem 3.2]{Vertman09}, summarizes the classical de Rham theorem on manifolds with boundary, 
cf. \cite[Remark after Proposition 4.2]{RS} and \cite[Theorem 4.1]{BL1}; strong ellipticity of the 
corresponding Laplace operators follows from \cite[Lemma 1.11.1]{Gi}. 

\begin{Thm}\label{Thm41} 
The Hilbert complexes $(\domr, \Dr)$ and $(\doma, \Da)$ are Fredholm and
the associated Laplacians $\triangle_{\textup{rel}}$ and $\triangle_{\textup{abs}}$ 
are strongly elliptic. The cohomologies $H^*(M,\partial M, E)$ and $H^*(M,E)$ of the Fredholm complexes 
$(\domr, \Dr)$ and $(\doma, \Da)$, respectively, can be computed from the following smooth subcomplexes, 
\begin{align*}
(\Omega^*_{\min}(M,E), \D), \quad &\Omega_{\min}^*(M,E):=\{\w \in \Omega^*(M,E)|\iota^*(\w)=0\}, \\
(\Omega^*_{\max}(M,E), \D), \quad &\Omega_{\max}^*(M,E):=\Omega^*(M,E),
\end{align*}
respectively, where $\iota: \partial M \hookrightarrow M$ 
denotes the natural inclusion of the boundary. Corresponding statement 
holds also for the complexes associated to the dual connection $\D'$. 
\end{Thm}

\subsection{The chirality operator} 
The Riemannian metric $g$ and a fixed orientation on $M$ define 
the Hodge star operator $*$ and the chirality operator ($r:= (m+1)/2$)
\begin{align}
\G :=i^r(-1)^{\frac{k(k+1)}{2}}*:\Omega^k(M,E)\to \Omega^{m-k}(M,E).
\end{align}
This operator extends to a self-adjoint involution on $L^2_*(M,E)$. 
The following properties of $\G$ are essential for the construction below, 
cf. \cite{Vertman09}.
 
\begin{Lem}\label{G-Lemma} 
The self-adjoint involution $\G$ on $L^2_*(M,E)$ maps $\dom(\Dr)$ to $\dom (\Da'^*)$,
and $\dom(\Da)$ to $\dom (\Dr'^*)$. With $\G$ restricted to appropriate domains, we have 
\begin{align*}
\G\Dr \G=\Da'^*, \quad 
\G\Da \G=\Dr'^*.
\end{align*}
\end{Lem}

\begin{Def}\label{domain}
We introduce the doubled Hilbert complexes
\begin{align*}
&(\domm, \DD):=(\domr, \Dr)\oplus (\doma, \Da),\\
&(\domm', \DD'):=(\domr', \Dr')\oplus (\doma', \Da').
\end{align*} 
\end{Def}

Similar to \eqref{equality}, we have the equality of domains
$$
\dom(\DD)=\dom(\DD'), \quad  \dom(\DD^*)=\dom (\DD'^*).
$$
The self-adjoint involution $\G$ gives rise to the "chirality operator"
\begin{align}\label{chirality}
\GG :=\left(\begin{array}{rr} 0 & \G \\ \G & 0 \end{array}\right)
\ \textup{on} \ L^2_*(M,E)\oplus L^2_*(M,E).
\end{align}
An immediate consequence of Lemma \ref{G-Lemma} is the following
\begin{Prop}
The chirality operator $\GG$ acts as
\begin{align*}
\GG|_{\dom(\DD)}: \dom(\DD) \to \dom(\DD^*), \quad
\GG|_{\dom(\DD^*)}: \dom(\DD^*) \to \dom(\DD).
\end{align*}
Moreover we have the relation $\GG \DD=\DD'^* \GG$.
\end{Prop}

\subsection{The odd signature operator}
We now apply the concepts of \cite{BraKap:refined} to our new setup
and define the odd-signature operator of the Hilbert complex 
$(\domm, \DD)$ by 
\begin{align}\label{odd-signature} 
\B:=\GG\DD+\DD\GG, \quad \dom(\B)=\dom(\DD)\cap \dom (\DD^*).
\end{align}

By \cite{Vertman09} the odd signature operator $\B$ is strongly elliptic 
with discrete spectrum and an Agmon angle $\theta \in (-\pi, 0)$.

\subsection{Spectral decomposition}
Consider for any $\lambda \geq 0$ the spectral projection of $\B^2$ onto 
eigenspaces with eigenvalues of absolute value in the interval $[0,\lambda]$:
$$
\Pi_{\B^2, [0,\lambda]}:=\frac{i}{2\pi}\int_{\gamma(\lambda)}(\B^2-x)^{-1}dx,
$$
with $\gamma(\lambda)$ being a closed counterclockwise circle around the origin surrounding 
eigenvalues of absolute value in $[0,\lambda]$.
By the analytic Fredholm theorem, the range of the projection 
lies in $\dom (\B^2)$ and the projection commutes with $\B^2$.
Moreover, $\Pi_{\B^2, [0,\lambda]}$ is of finite rank and the decomposition
\begin{align}\label{decomp-L-2}
L^2_*(M,E\oplus E)=\textup{Image}\,\Pi_{\B^2, [0,\lambda]}\oplus \textup{Image}\,(\one - \Pi_{\B^2, [0,\lambda]}),
\end{align}
is a direct sum decomposition into closed subspaces of the Hilbert space $L^2_*(M,E\oplus E)$.
Note that if $\alp$ is unitary and hence $\B^2$ is self-adjoint, 
the projection $\Pi_{\B^2,[0,\lambda]}$ is orthogonal. 
\eqref{decomp-L-2} induces a decomposition of $\domm$
$$
\domm=\domm_{[0,\lambda]}\oplus \domm_{(\lambda, \infty)}.
$$ 
Since $\DD$ commutes with $\B, \B^2$ and hence also with $\Pi_{\B^2, [0,\lambda]}$, 
we obtain a decomposition of $(\domm, \DD)$ into subcomplexes
\begin{equation}\label{decomposition}
\begin{split}
&(\domm, \DD)=(\domm_{[0,\lambda]}, \DD_{[0,\lambda]})\oplus 
(\domm_{(\lambda, \infty)}, \DD_{(\lambda, \infty)}) \\ 
&\textup{where} \ \DD_{\mathcal{I}}:=\DD|_{\domm_{\mathcal{I}}} 
\ \textup{for} \ \mathcal{I}=[0,\lambda] \ \textup{or} \ (\lambda, \infty).
\end{split}
\end{equation}
The chirality operator $\GG$ commutes with $\B, \B^2$ and respects the 
decomposition \eqref{decomposition} so that 
\begin{equation}\label{66}
\GG = \GG_{[0,\lambda]} \oplus \GG_{(\lambda, \infty)}, \quad \B=\B^{[0,\lambda]}\oplus \B^{(\lambda, \infty)}.
\end{equation}

\begin{Prop}\label{bijective}\label{cohomology} \cite[Corollary 3.14 and 3.15]{Vertman09}. 
The operator $\B^{(\lambda, \infty)}$, $\lambda \geq 0$ is bijective.
The complex $(\domm_{(\lambda, \infty)}, \DD_{(\lambda, \infty)})$ is acyclic and 
$$
H^*(\domm_{[0,\lambda]}, \DD_{[0,\lambda]})\cong H^*(\domm, \DD).
$$
\end{Prop}

\subsection{The refined torsion element}
Recall the notion of a determinant lines of a finite dimensional complex $(C^*,\partial_*)$ and of its cohomology. Set
\[
	\begin{aligned}
		\textup{Det}\, C^*\ &=\ 
		\bigotimes\limits_k \det\, (C^k)^{(-1)^k},\\
		\textup{Det}\,H^*(C^*,\partial_*)\ &=\ 
		\bigotimes\limits_k \det H^k(C^*,\partial_*)^{(-1)^k}, 
	\end{aligned}
\]
where for a vector space $V$ we denote by $\det V$ its top exterior power  and 
the $(-1)$ upper index denotes the dual vector space.  
We follow \cite[Section 1.1]{BraKap:refined-det} and define the canonical isomorphism 
\[
	\phi:\, \textup{Det}\, C^*\ \to\ \textup{Det}H^*(C^*,\partial_*)
\]
and the \emph{refined torsion element} 
of the complex $(\domm_{[0,\lambda]}, \DD_{[0,\lambda]})$
\begin{equation}\label{finite-torsion}
\begin{split}
\rho_{[0,\lambda]}\ :=\ \phi\big(\,c_0\otimes (c_1)^{-1}\otimes \cdots 
\otimes (c_r)^{(-1)^r} \otimes (\GG_{[0,\lambda]}c_r)^{(-1)^{r+1}}\otimes \cdots& 
\\ \cdots \otimes (\GG_{[0,\lambda]}c_1) \otimes (\GG_{[0,\lambda]}c_0)^{(-1)}\,\big)
\in \textup{Det}(H^*(\domm_{[0,\lambda]}, \DD_{[0,\lambda]})), &
\end{split}
\end{equation}
where $c_k\in \domm_{[0,\lambda]}$ 
are arbitrary elements of the determinant lines, we denote the extension of $\GG_{[0,\lambda]}$ 
to a mapping on determinant lines by the same letter, and for any $v\in \det\domm_{[0,\lambda]}$ 
the dual $v^{-1}\in \det (\domm_{[0,\lambda]})^{-1}\equiv 
\det (\domm_{[0,\lambda]})^*$ is the unique element such that $v^{-1}(v)=1$.

By Proposition \ref{cohomology} we can view $\rho_{[0,\lambda]}$ 
canonically as an element of $\textup{Det}(H^*(\domm, \DD))$, which we do henceforth.

\subsection{The graded determinant}
The fundamental part of the construction is the graded determinant. 
The operator $\B^{(\lambda, \infty)},\lambda \geq 0$ is bijective by Proposition \ref{bijective} 
and hence by injectivity (put $\mathcal{I}=(\lambda, \infty)$ to simplify the notation) 
\begin{align}\label{kern}
\textup{ker}(\DD_{\mathcal{I}}\GG_{\mathcal{I}})\cap\textup{ker}(\GG_{\mathcal{I}}\DD_{\mathcal{I}})=\{0\}.
\end{align}
Moreover the complex $(\domm_{\mathcal{I}}, \DD_{\mathcal{I}})$ is acyclic by 
Proposition \ref{cohomology} and due to $\GG_{\mathcal{I}}$ being an involution on $\textup{Im}(1-\Pi_{\B^2,[0,\lambda]})$ we have

\begin{equation}\label{image1}
\begin{split}
\textup{ker}(\DD_{\mathcal{I}}\GG_{\mathcal{I}})&=\GG_{\mathcal{I}}
\textup{ker}(\DD_{\mathcal{I}})=\GG_{\mathcal{I}}\textup{Im}(\DD_{\mathcal{I}})=\textup{Im}
(\GG_{\mathcal{I}}\DD_{\mathcal{I}}), \\  
\textup{ker}(\GG_{\mathcal{I}}\DD_{\mathcal{I}})&=\textup{ker}(\DD_{\mathcal{I}})=
\textup{Im}(\DD_{\mathcal{I}})=\textup{Im}(\DD_{\mathcal{I}}\GG_{\mathcal{I}}).
\end{split}
\end{equation}

We have $\textup{Im}(\GG_{\mathcal{I}}\DD_{\mathcal{I}})+\textup{Im}(\DD_{\mathcal{I}}
\GG_{\mathcal{I}})=\textup{Im}(\B^{\mathcal{I}})$ and by surjectivity of $\B^{\mathcal{I}}$ 
we obtain from the last three relations above
\begin{align}\label{hilbert-decomposition}
\textup{Im}(1-\Pi_{\B^2, [0,\lambda]})=\textup{ker}(\DD_{\mathcal{I}}\GG_{\mathcal{I}})
\oplus\textup{ker}(\GG_{\mathcal{I}}\DD_{\mathcal{I}}).
\end{align}
Note that $\B$ leaves $\ker (\DD\GG)$ and $\ker (\GG\DD)$ invariant. Hence, we put
\begin{align*}
	\B^{+,(\lambda,\infty)}_{\textup{even}}\ :=\ 
	\B^{(\lambda,\infty)}\restriction \domm^{\textup{even}}\cap \ker (\DD\GG), 
	\\
	\B^{-,(\lambda,\infty)}_{\textup{even}}\ :=\
	\B^{(\lambda,\infty)}\restriction \domm^{\textup{even}}\cap \ker (\GG\DD).
\end{align*}
We arrive at a direct sum decomposition 
\[
\B^{(\lambda,\infty)}_{\textup{even}}\ =\ \B^{+,(\lambda,\infty)}_{\textup{even}} 
\oplus \B^{-,(\lambda,\infty)}_{\textup{even}}.
\]
By \cite{Vertman09}, there exists an Agmon angle $\theta\in (-\pi, 0)$ for $\B$, 
which is clearly an Agmon angle for the restrictions above, as well. 
For strongly elliptic boundary value problems $(D,B)$ of order $\w$ on $M$ with an Agmon angle $\theta \in (-\pi, 0)$,
the associated zeta-function is defined by
\begin{align*}
&\zeta_{\theta}(s, D_B):=\sum\limits_{\lambda \in \textup{Spec}(D_B)
\backslash \{0\}}m(\lambda)\cdot \lambda_{\theta}^{-s}, \quad \textup{Re}(s) > \frac{\dim M}{\w},
\end{align*}
where $\lambda_{\theta}^{-s}:=\textup{exp}(-s\cdot \log_{\theta}\lambda)$ and 
$m(\lambda)$ denotes the multiplicity of the eigenvalue $\lambda$. The zeta function
is holomorphic for $\textup{Re}(s) > \dim M / \w$ and admits a meromorphic 
extension to $\C$ with $s=0$ being a regular point. 
Consequently, the \emph{graded zeta-function} 
$$
\zeta_{gr,\theta}(s,\B^{(\lambda,\infty)}_{\textup{even}}):=
\zeta_{\theta}(s,\B^{+,(\lambda,\infty)}_{\textup{even}})-
\zeta_{\theta}(s,-\B^{-,(\lambda,\infty)}_{\textup{even}}), \ Re(s)\gg 0,
$$
is regular at $s=0$ and we may introduce the following
\begin{Def}\label{graded-determinant}
Let $\theta \in (-\pi, 0)$ be an Agmon angle for $\B^{(\lambda, \infty)}$. 
Then the graded determinant associated to $\B^{(\lambda, \infty)}$ and its 
Agmon angle $\theta$ is defined as follows: 
$$
\Det'\nolimits_{gr,\theta}(\B^{(\lambda, \infty)}_{\textup{even}}):= 
\textup{exp}\left(-\left.\frac{d}{ds}\right|_{s=0} \zeta_{gr,\theta}\left(s,\B^{(\lambda,\infty)}_{\textup{even}}\right)\right).
$$
\end{Def}

\subsection{Refined analytic torsion}

\begin{Prop}\label{rho-element} \cite{BraKap:refined, Vertman09}. 
The element $$\rho(\D, g):=\Det'\nolimits_{gr,\theta}(\B^{(\lambda, \infty)}_{\textup{even}})
\cdot \rho_{[0,\lambda]}\in \textup{Det}(H^*(\domm, \DD))$$ 
is independent of the choice of $\lambda \geq 0$ and choice of 
Agmon angle $\theta \in (-\pi, 0)$ for the odd-signature operator $\B^{(\lambda, \infty)}$. 
\end{Prop}

The construction of $\rho(\D, g)$ is in fact independent of the 
choice of a Hermitian metric $h^E$. Indeed, a variation of $h^E$ does not change 
the odd-signature operator $\B$ as a differential operator and 
different Hermitian metrics give rise to equivalent $L^2-$norms 
over compact manifolds. Hence $\dom (\B)$ is indeed 
independent of the particular choice of $h^E$. Independence of the choice of a 
Hermitian metric $h^E$ is essential, since for non-unitary flat vector bundles 
there is no canonical choice of $h^E$ and a Hermitian metric is fixed arbitrarily.

The refined analytic torsion is then obtained by studying the dependence of 
$\rho(\D, g)$ on the Riemannian metric. We cite the final result from \cite{Vertman09}.

\begin{Thm}\label{RAT-sign}
Let $(M,g)$ be an odd-dimensional oriented compact Riemannian manifold with boundary. 
Let $(E, \D, h^E)$ be a flat complex vector bundle over $M$. 
Consider the trivial vector bundle $M\times \C$ with a trivial connection $d$ and 
let $B:=\B(d)$ denote the associated odd-signature operator. 
$\eta(B)$ denotes the eta invariant of the even part $B_{\textup{even}}$. Put
\begin{align*}
\widehat{\xi}(d, g):=\frac{1}{2}\sum_{k=0}^m(-1)^{k}
\cdot k\cdot \zeta_{2\theta}(s=0,B^2\restriction \domm^k).
\end{align*}
Then the refined analytic torsion of $(M,E,\D)$
\begin{align}\label{rat-definition}
\rho_{\textup{an}}(\D):=\rho(\D, g)\cdot \exp\left[i\pi \, 
\textup{rk}(E)(\eta (B) + \widehat{\xi}(d, g))\right]
\end{align}
is modulo sign independent of the choice of $g$ in the interior of $M$.
\end{Thm}

\section{Holomorphic structure on the determinant line bundle}
\label{S:detlinebundle}

In the next step we interpret $\rho(\alp)$ as an analytic 
section of the determinant line bundle over the representation space. This requires a separate 
discussion of the analyticity for the refined torsion element 
and the graded determinant. The present section studies analyticity 
of the refined torsion element, while the next deals with analyticity of 
the graded determinant.

\subsection{The determinant line bundle}\label{SS:detlinebundle}
The space $\mathscr{R}:=\Rep$ of complex $n$-dimensional representations of $\pi_1=\p$ has a natural 
structure of a complex analytic space, cf., for example, \cite[\S13.6]{BraKap:refined} . 
For each $\alp\in \Rep$ we denote by $E_\alp$ the flat vector bundle over $M$ whose 
monodromy is equal to $\alp$. Then the disjoint union
\begin{equation}\label{E:Detbundle}
    \Dbundle \ := \ \bigsqcup_{\alp\in \mathscr{R}}\, \Det\big(H^\b(M,E_\alp)\big)
\otimes  \Det\big(H^\b(M,\d M, E_\alp)\big)
\end{equation}
has a natural structure of a holomorphic line bundle over $\mathscr{R}$, called the 
{\em determinant line bundle}. In this section we describe this
structure, using a CW-decomposition of $M$. Then, we show that the refined 
analytic torsion is  a nowhere vanishing holomorphic
section of $\Dbundle$. We continue in the notation fixed in \S \ref{explicit-unitary}.

\subsection{The combinatorial cochain complex}\label{SS:CKalp}
Fix a CW-decomposition $K=\{e_1\nek e_N\}$ of $M$. Let $\widetilde{K}$ denote 
the universal cover of $K$. Then the fundamental group $\pi_1(M)$ acts on $C^\b(\widetilde{K}, \C)$
from the right and $\C^n$ is a left module over the group ring $\C[\pi_1]$ via the 
representation $\alp$. Then the cochain complex $C^\b(K,\alp)$ is defined as
\begin{align}
C^\b(K,\alp) := C^\b(\widetilde{K}, \C) \otimes_{\C[\pi_1]} \C^n.
\end{align}

For each cell $e_j$, fix a lift $\tile_j$, a cell of the CW-decomposition of $\tilM$,
such that $\pi(\tile_j)= e_j$. By definition, the pull-back of the bundle 
$E_\alp$ to $\tilM$ is the trivial bundle $\tilM\times\CC^n\to
\tilM$. Hence, the choice of the cells $\tile_1\nek \tile_N$ identifies the 
cochain complex $C^\b(K,\alp)$ of the CW-complex $K$ with
coefficients in $\Ea$ with the complex
\begin{equation}\label{E:C(K,alp)}
    \begin{CD}
        0 \ \to \CC^{n\cdot k_0} @>{\pa_0(\alp)}>> 
\CC^{n\cdot k_1} @>{\pa_1(\alp)}>>\cdots @>{\pa_{m-1}(\alp)}>> 
\CC^{n\cdot k_m} \ \to \ 0,
    \end{CD}
\end{equation}
where $k_j\in \ZZ_{\ge0}$  ($j=0\nek m=\dim M$) is equal to the number of 
$j$-dimensional cells of $K$ and the differentials $\pa_j(\alp)$ are
$(nk_{j}\times{}nk_{j-1})$-matrices depending analytically on $\alp\in \Rep$.

The cohomology of the complex \eqref{E:C(K,alp)} is canonically isomorphic to $H^\b(M,\Ea)$. Let
\begin{equation}\label{E:phialp}
    \phi_{C^\b(K,\alp)}:\, \Det\big(\,C^\b(K,\alp)\,\big) \ {\longrightarrow} \ \Det\big(\,H^\b(M,\Ea)\,\big)
\end{equation}
denote the canonical isomorphism, cf. formula (2.13) of \cite{BraKap:refined-det} .

\subsection{A non-zero element of $\Det\big(H^\b(M,E_\alp)\big)$}\label{SS:sectionDetK}
The standard bases of $\CC^{n\cdot k_j}$ ($j=0\nek m$) define an element $c\in \Det\big(C^\b(K,\alp)\big)$, 
and, hence, an isomorphism $\psi_\alp:\,\CC \ {\longrightarrow} \ \Det\big(C^\b(K,\alp)\big)$ 
with $\psi_\alp(z) = z\cdot c$. Then for each $\alp\in \Rep$ we define
\begin{equation}\label{E:combsection}
    \sig(\alp) \ = \ \phi_{C^\b(K,\alp)}\big(\,\psi_\alp(1)\,\big) \ \in \ \Det\big(\,H^\b(M,\Ea)\,\big),
\end{equation}
a non-zero element of $\Det\big(\,H^\b(M,\Ea)\,\big)$. 
Of course, this element depends on the choice of the lifts $\tile_1\nek \tile_N$.

\subsection{A non-zero element of $\Det\big(H^\b(M,\d M,E_\alp)\big)$}\label{SS:sectionDetK'}
Let now $K'$ denote the CW-decomposition of $\d{}M$ induced by $K$. Then $K'\subset K$ 
and  the choice of the lifts $\tile_i$ made above identify the cochain complex $C^\b(K',\alp)$ 
of the CW-complex $K'$ with coefficients in $\Ea$ with the complex
\begin{equation}\label{E:C(K',alp)}
    \begin{CD}
        0 \ \to \CC^{n\cdot l_0} @>{\pa_0(\alp)}>> \CC^{n\cdot l_1} @>{\pa_1(\alp)}>>
\cdots @>{\pa_{m-1}(\alp)}>> \CC^{n\cdot l_{m-1}} \ \to \ 0,
    \end{CD}
\end{equation}
As above, the standard bases of $\CC^{n\cdot l_j}$ ($j=1,\ldots, m-1$) 
defines a canonical element of $\Det(C^\b(K',\alp))$ and an isomorphism 
$\psi'_\alp$ from $\CC$ onto $\Det\big(C^\b(K,\alp)\big)$. Thus we define 
\begin{equation}\label{E:combsection'}
    \sig'(\alp) \ = \ \phi_{C^\b(K',\alp)}\big(\,\psi'_\alp(1)\,\big) \ \in 
\ \Det\big(\,H^\b(\d M,\Ea)\,\big),
\end{equation}
 where $\alp'$ is the restriction of the representation $\alp$ to $\pi_1(\d{}M)$ 
and $E_\alp'$ is the restriction of $E_\alp$ to $\d{}M$.  
Consider the quotient complex 
\[
	C^\b(K,K',\alp) \ := \ C^\b(K,\alp)/C^\b(K',\alp).
\]
Using \eqref{E:C(K,alp)} and \eqref{E:C(K',alp)} we can identify 
$\Det\big(C^\b(K,K',\alp)\big)$ with $\CC$ thus constructing a map
\begin{equation}\label{E:phi''(c)}
    \psi''_\alp:\,\CC \ {\longrightarrow} \ \Det\big(C^\b(K,K',\alp)\big).
\end{equation}   

The cohomology of the complex $C^\b(K,K',\alp) $ is canonically 
isomorphic to the relative cohomology $H^\b(M,\d M,E_\alp)$. 
For each $\alp$ we define a non-zero element of 
$\Det\big(\,H^\b(M,\d M,\Ea)\,\big)$ by the formula 
\begin{equation}\label{E:combsection''}
    \sig''(\alp)\ = \   \phi_{C^\b(K,K',\alp)}\big(\,\psi''_\alp(1)\,\big) 
\ \in \ \Det\big(\,H^\b(M,\d M,\Ea)\,\big),
\end{equation}

\subsection{The holomorphic structure on $\Dbundle$}\label{SS:holonDet}
Recall that the elements $\sig(\alp)$ and $\sig''(\alp)$ are defined in 
\eqref{E:combsection} and \eqref{E:combsection''}. Consider the map
\begin{equation}\label{E:combsection_tau}
\begin{split}
    \tau:\, \alp \ &\mapsto \ \sig(\alp)\otimes\sig''(\alp)
    \\ &\in \ \Det\big(\,H^\b(M,\Ea)\,\big)\otimes \Det\big(\,H^\b(M,\d M,\Ea)\,\big),
    \end{split}
\end{equation}
where $\alp\in \Rep$, is a nowhere vanishing section of the 
determinant line bundle $\Dbundle$ over $\Rep$. 
\begin{Def}\label{D:holsection}
We say that a section $s(\alp)$ of\/ $\Dbundle$ is {\em holomorphic} 
if there exists a holomorphic function $f(\alp)$ on $\Rep$, such that
$s(\alp)= f(\alp)\cdot{}\tau(\alp)$.
\end{Def}
This defines a holomorphic structure on $\Dbundle$, which is independent 
of the choice of the lifts $\tile_1\nek \tile_N$ of $e_1\nek e_N$,
since for a different choice of lifts the section $\tau(\alp)$ will be multiplied 
by a constant. In the next subsection we show that this
holomorphic structure is also independent of the CW-decomposition 
$K$ of $M$.

\subsection{The Farber-Turaev torsion}\label{SS:Turaevtorsion}
The choice of the lifts $\tile_1\nek \tile_N$ of $e_1\nek e_N$ determines an 
{\em Euler structure} $\eps$ on $M$, while the ordering of the cells
$e_1\nek e_N$ determines a cohomological orientation $\gro$, cf. \cite[\S20]{Turaev01}. 
Moreover, every Euler structure and every cohomological
orientation can be obtained in this way. The Farber-Turaev torsion $\rtt(\alp)$, 
corresponding to the pair $(\eps,\gro)$, is, by definition,
\cite[\S6]{FarberTuraev00}, equal to the element $\sig(\alp)$ defined in \eqref{E:combsection}. 
Since the Farber-Turaev torsion is independent of the choice of the CW-decomposition of $M$, 
cf.  \cite{Turaev90,FarberTuraev00},  we conclude that the element $\sig(\alp)$ is also 
independent of the CW-decomposition, but only depends on the Euler structure 
and the cohomological orientation. 

The lifts $\tile_j$ and the ordering of the cells also defines an Euler structure 
$\eps'$ and a cohomological orientation $\gro'$ of $\d{}M$. The element 
$\sig'(\alp)$ defined in \eqref{E:combsection'} is equal to the Farber-Turaev 
torsion $\rho_{\eps',\gro'}(\alp')$ where $\alp'$ is the restriction of the 
representation $\alp$ to $\pi_1(\d{}M)$.  Let
\begin{equation}\label{E:fusionexact}
	\mu:\, \Det\big(C^\b(K',\alp)\big)\otimes\Det\big(C^\b(K,K',\alp)\big)\ 
\longrightarrow \ \Det\big(C^\b(K,\alp)\big)
\end{equation}
denote the fusion isomorphism, cf. \cite[\S2.6]{BraKap:refined-det} . 
Using this isomorphism we define the map
\begin{equation}\label{E:fusioncoh}
 \begin{aligned}
 \nu:\,  \Det\big(H^\b(\d M,&\Ea)\big)\otimes\Det\big(H^\b(M,\d M,\Ea)\big)\ 
\longrightarrow \ \Det\big(H^\b(M,\Ea)\big),\\
 \nu\ &:= \ \phi_{C^\b(K,\alp)}\circ\mu\circ \big(\,\phi_{C^\b(K',\alp)}^{-1}
\otimes \phi_{C^\b(K,K',\alp)}^{-1}\,\big).
 \end{aligned}
\end{equation}
It follows from the construction that the elements $\sig(\alp)$, $\sig'(\alp)$ 
and $\sig''(\alp)$ defined in \eqref{E:fusionexact}, \eqref{E:combsection}, 
\eqref{E:combsection'}, and  \eqref{E:combsection''}  satisfy the equality
\begin{equation}\label{E:sigsigsig}
	\nu\big(\,\sig'(\alp)\otimes \sig''(\alp)\,\big) \ = \ \pm \, \sig(\alp),
\end{equation}
where the sign depends only on the dimensions of the spaces $C^\b(M,\alp)$ 
and $C^\b(\d{}M,\alp)$ but not on the representation $\alp$. Since $\sig(\alp)$ 
and $\sig'(\alp)$ are independent of the CW-decomposition, it follows that 
$\sig''(\alp)$ is independent of CW-decompositon up to a sign. In fact $\sigma''(\alp)$ 
can be considered as the definition of the relative Farber-Turaev torsion. It follows 
that the section $\tau(\alp)$ defined in \eqref{E:combsection_tau} is also 
an independent of the CW-decomposition up to a sign. Hence so is  the 
holomorphic structure defined in \refd{holsection}.

\subsection{The acyclic case}\label{SS:Turaevacyclic}
Let $\Repo\subset\Rep$ denote the set of representations such that 
$H^\b(M,E_\alp)=0$, $H^\b(\d M,E_{\alp'})=0$, and $H^\b(M,\d M,E_\alp)=0$. 
Then the determinant lines $\Det(H^\b(M,E_\alp))$, $\Det(H^\b(M,\d M,E_{\alp'}))$, 
and $\Det(H^\b(M,\d M.E_\alp))$ are canonically isomorphic to $\CC$. Hence, the 
Farber-Turaev torsions $\rho_{\eps,\gro}(\alp)$ and  $\rho_{\eps',\gro'}(\alp)$ can 
be viewed as a complex-valued functions on $\Repo$. It is easy to see, cf. 
\cite[Theorem 4.3]{BurgheleaHaller_Euler}, that these functions are holomorphic on 
$\Repo$. Moreover, they are rational functions on $\Rep$, all whose poles are in 
\[
	\Rep\backslash\Repo.
\]
In particular, the holomorphic structure on $\Dbundle$, which we defined above, 
coincides, when restricted to $\Repo$, with the natural holomorphic structure 
obtained from the canonical isomorphism
\[
	\Dbundle|_{\Repo}\ \simeq\ \Repo\times\CC.
\]

\section{The graded determinant as a holomorphic function}
\label{S:drdetholom}

Fix $\alp_0\in \Rep$. Fix a number $\lam\ge0$ which is not in the spectrum of the square 
$\calB_{\alp_0}^2$ of the odd signature operator $\calB_{\alp_0}=\calB(\n_{\alp_0},g)$. 
Then there is a neighborhood $U_\lam\subset \Rep$ of $\alp_0$ such that $\lam$ is not an 
eigenvalue of $\calB_\alp^2$ for all $\alp\in U_\lam$. Denote by $\calB_\alp^{(\lam,\infty)}$ 
the restriction of $\calB_\alp$ to the spectral subspace of $\calB_\alp^2$ corresponding to 
the spectral set $(\lam,\infty)$. Then $\calB_\alp^{(\lam,\infty)}$ is an invertible operator. 
Let $\tet\in (-\pi/2,0)$ be an Agmon angle for $\calB_{\alp_0}^{(0,\lam)}$ and assume that there 
are no eigenvalues of $\calB_{\alp_0}^{(\lam,\infty)}$  in the solid angles $L_{(-\pi/2,\lam]}$ 
and $L_{(\pi/2,\lam+\pi/2]}$. Then there exists a neighborhood $U_{\lam,\tet}\subset U_\lam$ 
of $\alp_0$ such that $\tet$ is also an Agmon angle for $B^{(\lam,\infty)}_\alp$ for all 
$\alp\in U_{\lam,\tet}$.  In this section we prove that   the graded determinant 
$\Detgrtet(\calB^{(\lambda, \infty)}_{\alp,\even})$   is  a holomorphic function on $U_{\lam,\tet}$.
Our main result in this section is the following
\begin{Thm}\label{T:analytic}
Let $\calO\subset \CC$ be a connected open neighborhood of $0$. Let
\[
	\gam:\calO\ \to\  U_{\lam,\tet}\subset \Rep
\] 
be a holomorphic curve such that  $\gam(0)=\alp_0$.  Then the function
\begin{equation}\label{E:weakanalytic}
   z \ \mapsto \Detgrtet\big(\,\calB_{\gam(z),\even}^{(\lam,\infty)}\,\big)
\end{equation}
is holomorphic in a neighborhood of $0$.
\end{Thm}

An immediate consequence is the following

\begin{Cor}\label{C:analytic}
Suppose  $V\subset U_{\lam,\tet}$ is an open subset  such that all points 
$\alp\in V$ are regular points of the complex algebraic set $\Rep$. Then the map
\begin{equation}\notag
    \Det:\, V \ \longrightarrow \ \CC, \qquad \Det:\,\alp \  \mapsto\ \Det(\alp) \ 
:= \ \Detgrtet(\B_{\alp,\even}^{(\lam,\infty)}).
\end{equation}
is holomorphic.
\end{Cor}

\begin{proof}
By Hartogs' theorem (cf., for example, \cite[Th.~2.2.8]{HormanderSCV}), a 
function on a smooth algebraic variety is holomorphic if its restriction to
each holomorphic curve is holomorphic. Hence, the corollary follows immediately from \reft{analytic}.
\end{proof}

The rest of this section is occupied with the proof of Theorem \ref{T:analytic}.

\subsection{A germ of connections}\label{SS:germconnections}
Let us introduce some additional notations. Let $E$ be a vector bundle over 
$M$ and let $\calC(E)$ denote the affine space of (not necessarily flat) connections on $E$. 
We endow  $\calC(E)$ with the  the Fr\'echet topology on $\calC(E)$ introduced in 
Section~13.1 of \cite{BraKap:refined} . 

Fix a base point $x_*\in M$ and let $E_{x_*}$ denote the fiber of $E$ over $x_*$. 
We will identify $E_{x_*}$ with $\CC^n$ and $\pi_1(M,x_*)$ with $\pi_1(M)$.

For $\n\in \calC(E)$ and a closed  path $\phi:[0,1]\to M$ with $\phi(0)=\phi(1)= x_*$, 
we denote by $\Mon_\n(\phi)\in \End{}E_{x_*}\simeq \Mat_{n\times n}(\CC)$ the 
monodromy of $\n$ along $\phi$. Note that, if $\n$ is flat then $\Mon_\n(\phi)$ 
depends only on the class $[\phi]$ of $\phi$ in
$\pi_1(M)$. Hence, if $\n$ is flat, then the map $\phi\mapsto \Mon_\n(\phi)$ 
defines an element of $\Rep$, called  the {\em monodromy representation} of $\n$.

Suppose now that $\calO\subset \CC$ is a connected open neighborhood of $0$. For simplicity we also assume that $\calO$ is convex. Let
\[
  \gam:\,\calO\ \to\ \Rep
\]
be a holomorphic curve with $\gam(0)=\alp_0$. The operator $\calB_{\gam(z)}$ 
is constructed using a flat connection $\n_{\gam(z)}$ whose monodromy  is equal 
to $\gam(z)$. Unfortunately, there is no a canonical choice of such connection. 
Though the graded determinant $\Detgrtet\big(\,\calB_{\gam(z),\even}^{(\lam,\infty)}\,\big)$ 
is independent of this choice, to study the dependence of this determinant on $z\in \calO$ 
we need to choose a family of connections $\n_{\gam(z)}$. The main difficulty in the 
proof of Theorem~\ref{T:analytic} is that  it is not clear whether there exists a 
\emph{holomorphic} family $\n_{\gam(z)}$ with $\Mon_{\n_{\gam(z)}}=\gam(z)$. 
We shall now explain how to circumvent this difficulty. 

By Proposition~4.5 of \cite{GoldmanMillson88}, all the bundles $E_{\gam(z)}, \ z\in \calO$, are isomorphic to each other. Moreover, we have the following lemma.


\begin{Lem}\label{L:diff connection}
There exists a vector bundle $E\to M$ and a real differentiable 
family of flat connections $\n_{\gam(z)},\ z\in \calO$, on $E$, such that the monodromy representation of $\n_{\gam(z)}$ is equal to $\gam(z)$ for all $z\in \calO$. 
\end{Lem}
\begin{proof}
By Lemma~3.3 of  \cite{BairdRamras15} there exists a smooth vector bundle 
\[
	\tilde{E}\ \to\ M\times\calO
\] 
and a smooth connection $\tilde{\nabla}$ on $\tilde{E}$, whose restriction $\tilde{\n}_z$ to 
\[
	\tilde{E}\big|_{M\times\{z\}}\ \to \ M\times\{z\}
\] 
is flat for each $z\in \calO$, and such that the monodromy of $\tilde{\n}_z$ is equal to  $\gamma(z)$. 

Set $E:= \tilde{E}\big|_{M\times\{0\}}$ and let 
\[
	\Phi_z:\, E\ \to \tilde{E}\big|_{M\times\{z\}}
\]
denote the parallel transport along the intervals 
\[
	\big\{\,(m,tz):\,m\in M, \ t\in [0,1]\,\big\}\ \subset\  M\times\calO.
\]
Then 
\[
	\nabla_{\gamma(z)}\ := \ \Phi_z^{-1}\circ \tilde{\nabla}_z\circ \Phi_z,
	\qquad z\in \calO,
\]
is a smooth family of connections on $E$, and the monodromy of $\nabla_{\gamma(z)}$ is equal to $\gamma(z)$. 
\end{proof}

Furthermore, Lemma~B.6 of
\cite{BraKap:refined}  shows that the family $\n_{\gam(z)}$ can be chosen so 
that there exists a one-form $\ome\in \Ome^1(M,\End{E})$ such that
\begin{equation}\label{E:n0+uome}
    \n_{\gam(z)} \ = \ \n_{\alp_0} \ + \  z\cdot \ome \ + \ o(z),
\end{equation}
where $o(z)$ is understood in the sense of the Fr\'echet topology on $\calC(E)$ 
introduced in Section~13.1 of \cite{BraKap:refined}, and always refers to the behavior as $z\to 0$.

\begin{Rem}
Lemma~\ref{L:diff connection} asserts that the family $\nabla_{\gamma(z)}$ can be chosen to be  {\em real differentiable} at every point $z\in \calO$. We don't know whether it can be chosen to be complex differentiable on the whole set $\calO$. However, the equation \eqref{E:n0+uome} implies that it can be chosen to be complex differentiable at 0. 
\end{Rem}

\subsection{The family of projections}\label{SS:familyproj}
Let $\Flat(E)\subset \calC(E)$ denote the set of flat connections on $E$ and consider a curve 
$\n(z)\in \Flat(E)$, $z\in \calO$, of connections such that $\Mon \n(0)=\alp_0$. 
We will assume that it is complex differentiable at 0 in the sense  that 
\begin{equation}\label{E:n(z)}
	\n(z) \ = \ \n_{\alp_0} \ + \ z\cdot \ome \ + \ o(z), \quad z\to 0,
\end{equation}
where $\ome\in \Ome^1(M,\End E)$. We denote by 
\[
	\DD(z) \ := \  \n(z)_{\operatorname{min}}\oplus \n(z)_{\operatorname{max}}
\] 
and by $\B(z)$ the corresponding odd signature operator the sense of \S \ref{explicit-unitary}.  
Let $P(z)=\Pi_{\B^2(z), (\lambda, \infty)}$ denote the spectral projection of the operator 
$\B(z)^2$ onto the subspace $\widetilde{\calD}_{(\lam,\infty)}$ spanned by eigenforms of 
$\B(z)^2$ with eigenvalues in $(\lam,\infty)$. 

From \eqref{hilbert-decomposition} we conclude that the space $L^2_*(M,E\oplus{}E)$ of $L^2$-differential forms is a direct sum
\[
	L^2_*(M,E\oplus{}E)\ =\ \IM\big(I-P(z)\big)\oplus\overline{\IM \big(P(z)\DD(z)\big)}\oplus\overline{\IM\big(P(z)\GG \DD(z)\big)}.
\]
Consider the corresponding orthogonal projections
\begin{equation}\label{E:Pi+-}\begin{aligned}
	P_+(z)&:L^2_*(M,E\oplus{}E)\ \longrightarrow \  \overline{\IM \big(P(z)\DD(z)\big)}, \\ 
	P_-(z)&:L^2_*(M,E\oplus{}E)\ \longrightarrow \  \overline{\IM\big(P(z)\GG \DD(z)\big)}.
 \end{aligned}
\end{equation}
 Then 
\[
	\IM P(z)\ = \ \IM P_+(z)\oplus\IM P_-(z).
\]

\lem{familyofprojections} \label{fam-proj}
There exist bounded operators $A_\pm$ on $L^2_*(M,E\oplus{}E)$ with 
\begin{equation}\label{E:familyofprojections}
\begin{split}
	P_\pm(z) \ = \ P_\pm(0)\ &+\ z\, \Big(\,P_\pm(0)A_\pm\big(\Id-P_\pm(0)\big)\\ &+ 
\ \big(\Id-P_\pm(0)\big)A_\pm P_\pm(0)\,\Big) + \ o(z).
\end{split}
\end{equation}
\elem
\prf
For small enough $z\in \calO$  the operators $P_\pm(z)$ depend smoothly on $z$, 
there exist bounded operators $A_\pm$ and $A$ such that 
\begin{equation}\label{E:familyofprojections1}
\begin{split}
	P_\pm(z) \ &= \ P_\pm(0)\ +\ z\,A_\pm\ + \ o(z), \\
      P(z) \ &= \ P(0)\ +\ z\,A\ + \ o(z).
\end{split}
\end{equation}
Using the decomposition 
\begin{multline}\label{E:decofA}
	A\ = \ P_\pm(0)AP_\pm(0)\ +\ \big(\Id-P_\pm(0)\big)AP_\pm(0)\\ +
\ P_\pm(0)A\big(\Id-P_\pm(0)\big)\ +\ 
	\big(\Id-P_\pm(0)\big)A\big(\Id-P_\pm(0)\big)
\end{multline} 
and the equality $P_\pm(z)^2=P_\pm(z)$ we obtain
\begin{multline}\label{E:familyofprojections2}
	P_\pm(z) \ = \ P_\pm(z)^2\ = \ P_\pm(0)\ +\ 
	2 z\,P_\pm(0)AP_\pm(0)\\ +\ z\,\big(\Id-P_\pm(0)\big)AP_\pm(0)\ +\ 
 z\,P_\pm(0)A\big(\Id-P_\pm(0)\big)\ + \ o(z).
\end{multline}

Comparing \eqref{E:familyofprojections1} and \eqref{E:familyofprojections2} 
and using \eqref{E:decofA} we conclude that 
\[
	P_\pm(0)AP_\pm(0) \ = \ \big(\Id-P_\pm(0)\big)A\big(\Id-P_\pm(0)\big) \ =\ 0.
\]
The equality \eqref{E:familyofprojections} follows now from 
\eqref{E:familyofprojections1} and \eqref{E:decofA}.
\eprf

\subsection{The partial derivatives of the graded determinant}\label{SS:partialDett}
In terms of the projections $P_\pm(z)$ introduced in the previous section the odd signature 
operator $\calB^{(\lam,\infty)}$ can be written in the form
\begin{equation}\label{E:BP(z)}
\begin{split}
	\calB^{(\lam,\infty)}(z)\ &= \ \big(\,\DD(z)\GG  \ + \ \GG\DD(z)\,\big)\, P(z)\\ =& \ 
	\DD(z)\GG P_+(z)  \ + \ \GG\DD(z) P_-(z).
	\end{split}
\end{equation}
Hence, we may write
\begin{equation*}
	\Detgrtet\big(\calB^{(\lam,\infty)}_\even(z)\big)\ = \ \frac{
	\Det'_\tet\big(\DD(z)\GG P_+(z)\restriction L_\even^2(M,E\oplus E)\big)}
	{\Det'_\tet\big(-\GG\DD(z) P_-(z)\restriction L_\even^2(M,E\oplus E)\big)}
\end{equation*}
Consider a curve $\kap:(-1,1)\to \calO$ such that $\kap(0)=\alp_0$. To simplify the notation, set 
\begin{equation}\label{E:Dpm}
\begin{split}
 	&D_+(z)\ := \ \DD(z)\GG P_+(z)\restriction L_\even^2(M,E\oplus E), \\
 	&D_-(z)\ := \ -\GG\DD(z)  P_-(z)\restriction L_\even^2(M,E\oplus E),
 		\end{split}
\end{equation}
and also 
\begin{equation}\label{E:tn'P'}
 	\DD'\ :=\ \frac{d}{dt} \DD(\kap(t))\big|_{t=0},\ P'_\pm\ :=\ \frac{d}{dt}  P_\pm(\kap(t))\big|_{t=0},
	\ D_\pm'\ = \ \frac{d}{dt} D_\pm\big(\kap(t)\big)\big|_{t=0}.
\end{equation}
With this notation we have 
\begin{equation}\label{E:detD=D+/D-}
 	\Detgrtet\big(\calB^{(\lam,\infty)}_\even(z)\big) \ = \ \frac{\Det_\tet' \big(D_+(z)\big)}{\Det_\tet'\big(D_-(z)\big)}.
\end{equation}
We consider the function
\begin{equation}\label{E:F(z)}
	F(z) \ =  \ \frac{
	\Det'_\tet\big(\DD(z)\GG P_+(\alp_0)\restriction L_\even^2(M,E\oplus E)\big)}
	{\Det'_\tet\big(-\GG\DD(z) P_-(\alp_0)\restriction L_\even^2(M,E\oplus E)\big)}.
\end{equation}
Notice that the right hand side of \eqref{E:F(z)} is similar to the right hand side of \eqref{E:detD=D+/D-} but $P_\pm(z)$ is replaced by $P_\pm(0)= P_\pm(\alp_0)$. In particular, \[
	F(0)\ =\  \Detgrtet\big(\calB^{(\lam,\infty)}_\even(0)\big).
\]

\lem{partialDet}
Then for any curve $\kap:(-1,1)\to \calO$ with $\kap(0)=\alp_0$ we have
\begin{equation}\label{E:partialDet}
	\frac{d}{dt}\Big|_{t=0}\log\Detgrtet\big(\calB^{(\lam,\infty)}_\even(\kap(t))\big)\ = \ \frac{d}{dt}\Big|_{t=0}\log F\big(\kap(t)\big).
\end{equation}
\elem
\prf
 Using Lemma \ref{fam-proj} and the equality
 \begin{equation}\label{E:NGP=PnG}
 	\DD(\alp_0)\GG \, P_+(\alp_0)\ =\ P_+(\alp_0)\,\DD(\alp_0)\GG,
\end{equation}
we obtain
 \begin{equation}\label{E:D+'=}
 \begin{split}
 	D_+'\ = \ \DD' \GG P_+(\alp_0) \ &+ \  P_+(\alp_0)\,\DD(\alp_0)\GG\, A_+\,\big(\Id-P_+(\alp_0)\big)\ \\ &+ \ 
	\big(\Id-P_+(\alp_0)\big)\,\DD(\alp_0)\GG \,A_+ P_+(\alp_0).
	\end{split}
 \end{equation}
By the variation formula for the logarithm of the determinant, cf., for example, Section~3.7 of \cite{BFK92}, we have 
\begin{multline}\label{E:ddtD+long}
	\frac{d}{dt}\Big|_{t=0}\log\Det'_\tet D_+(\kap(t)) \ = \ 
	\Tr D_+^{-s-1}(\alp_0) D_+'\Big|_{s=0} \\ = \
	\Tr D_+^{-s-1}(\alp_0) \DD' \GG P_+(0)\Big|_{s=0} \ + \ \Tr D_+^{-s-1}(\alp_0) \DD_{\alp_0} \GG P_+'\Big|_{s=0}.
\end{multline}
Using \eqref{E:familyofprojections} and the fact that  the operators $D_+(0)$ and $\DD_{\alp_0}\GG$ commute with $P_+(0)$ we conclude that 
\begin{multline}\label{E:secondvanish}
	\Tr D_+^{-s-1}(0) \DD_{\alp_0} \GG P_+' \ = \ \Tr \big(\Id-P_+(0)\big)\,D_+^{-s-1}(0) \DD_{\alp_0} \GG A_+P_+(0) \\ + \ 
	\Tr P_+(0)D_+^{-s-1}(0) \DD_{\alp_0} \GG A_+\,\big(\Id-P_+(0)\big) \ = \ 0.
\end{multline}
Combining \eqref{E:ddtD+long} and \eqref{E:secondvanish} we conclude that 
\begin{equation}\label{E:ddtD+long2}
\begin{split}
	&\frac{d}{dt}\Big|_{t=0}\log\Det'_\tet D_+(\kap(t)) \ = \ 	\Tr D_+^{-s-1}(\alp_0) \DD' \GG P_+(0)\Big|_{s=0}
	 \\ =\ &\frac{d}{dt}\Big|_{t=0}\Det'_\tet\big(\DD(\kap(t))\GG P_+(\alp_0)\restriction L_\even^2(M,E\oplus E)\big).
	 \end{split}
\end{equation}
Similarly, 
\begin{equation}\label{E:ddtD_-long}
\begin{split}
	&\frac{d}{dt}\Big|_{t=0}\log\Det'_\tet D_-(\kap(t)) \\ = \ &\frac{d}{dt}
	\Big|_{t=0}\log\Det'_\tet\big(-\GG\DD(\kap(t)) P_-(\alp_0)\restriction L_\even^2(M,E\oplus E)\big).
	\end{split}
\end{equation}
From \eqref{E:detD=D+/D-}, \eqref{E:F(z)}, \eqref{E:ddtD+long2} and \eqref{E:ddtD_-long} we obtain \eqref{E:partialDet}.
\eprf

\begin{Prop}\label{P:analyticn}
Let $\n(z)$ $(z\in \calO)$ be a family of flat connections such that as $z\to 0$
\begin{equation}\label{E:analyticn}
	\n(z)\ = \ \n_{\alp_0}\ +\ z\cdot\ome \ + \ o(z).
\end{equation}
Then the function 
\begin{equation}\label{E:weakanalytic1}
   z \ \mapsto f(z)\ := \ \Detgrtet\big(\,\calB_{\even}^{(\lam,\infty)}(\n_{\gam(z)},g)\,\big)
\end{equation} 
is complex differentiable at zero. In other words, there exists a complex number $A$ such that 
\begin{equation}\label{E:detiscomplexdif}
\begin{split}
	\Detgrtet\big(\,\calB_{\even}^{(\lam,\infty)}(\n_{\gam(z)},g)\,\big)\ &= \ \Detgrtet\big(\,\calB_{\even}^{(\lam,\infty)}(\n_{\gam(0)},g)\,\big)
	\\ &+\ A\cdot z \ + o(z).
	\end{split}
\end{equation}
\end{Prop}
\begin{proof}
By \eqref{E:partialDet} it is enough to show that 
\begin{equation}\label{E:detiscomplexdif2}
	F(z) \ = \ F(0) \ + A\cdot z \ + \ o(z).
\end{equation}

Set $z=x+iy$. Using the variation formula for the logarithm of the determinant as in the proof of Lemma~\ref{L:partialDet}, one easily sees that 
\[
	i\frac{\d}{\d x} \log F(z)\ = \ \frac{\d}{\d y}\log F(z),
\] 
which is equivalent to \eqref{E:detiscomplexdif2}.
\end{proof}

\subsection{Proof of Theorem~\ref{T:analytic}}\label{SS:pranalytic}
It follows from Proposition~\ref{P:analyticn} that the function $$z\mapsto \Detgrtet\big(\calB_{\even}^{(\lam,\infty)}(\n_{\gam(z)},g)\big),$$ is complex differentiable at 0. 

Let $a\in \calO$ be such that $\gam(a)\in U_{\lam,\tet}$. By making a change of variables $\zet=z-0$ we conclude that this function is also complex 
differentiable at $a$. Hence, this function is holomorphic in $\gam^{-1}(U_{\lam,\tet})$.
\hfill$\square$

\section{Refined analytic torsion as a holomorphic section}\label{S:holomorphic}

In this section we show that the refined analytic torsion $\rat$ is a non-vanishing holomorphic section of\/ $\Dbundle$. More
precisely, our main result is the following

\subsection{Weakly holomorphic section}\label{SS:weakly hol}
Recall from \cite[p.~148]{Gunning70-AnVarieties1} that a continuous function on a singular space $X$ is called {\em weakly holomorphic} if its restriction to the set of regular points of $X$ is holomorphic. Such functions have many properties of analytic functions on $X$. In particular, they are all meromorphic . We refer to \cite{Gunning70-AnVarieties1}, \cite[\S1e]{Gunning74-AnVarieties2} for the properties of the weakly holomorphic functions. 

We say that a section $s(\alp)$ of\/ $\Dbundle$ is {\em weakly holomorphic} 
if there exists a weakly holomorphic function $f(\alp)$ on $\Rep$, such that
$s(\alp)= f(\alp)\cdot{}\tau(\alp)$.

\begin{Thm}\label{T:holomsection}
The refined analytic torsion $\rat$ is a weakly holomorphic section of the determinant bundle $\Dbundle$.
In particular, the restriction of $\rat$ to the set $\Repo$ of acyclic representations, viewed as a complex-valued function via the canonical
isomorphism
\[
    \Dbundle|_{\Repo}\ \simeq\ \Repo\times\CC,
\]
is a weakly holomorphic function on $\Repo$.
\end{Thm}

\subsection{Reduction to a finite dimensional complex}\label{SS:redtofd}
Let $\alp_0\in\linebreak\Rep$. Fix a Riemannian metric $g$ on $M$ and a number $\lam\ge0$ such that there are no eigenvalues of\/ $\B(\n_{\alp_0},g)^2$ with absolute value equal to $\lam$. Let $\tet$ be an Agmon angle for $\B(\n_{\alp_0},g)^2$ and let $U_{\lam,\alp}\subset \Rep$ be as in \refs{drdetholom}. By Corollary~\ref{C:analytic} the function $\alp\mapsto \Detgrtetnp(\B_\even^{(\lam,\infty)}(\na,g))$ is weakly holomorphic  on $U_{\lam,\tet}$. It follows now from the definition of the refined analytic torsion that to prove \reft{holomsection} it is enough to show that 
\[
    \alp \mapsto \ \rho_{[0,\lam]}\ \equiv\ \rho_{\GG_{[0,\lam]}}(\na,g)
\]
is a weakly holomorphic section of $\Rep$.  By the definition of the holomorphic structure on the bundle $\Dbundle$, cf. \refd{holsection}, this means that the function 
\[
    \alp \mapsto \ \frac{\rho_{\GG_{[0,\lam]}}(\na,g)}{\tau(\alp)},
\]
is continuous at $\alp_0$ and is holomorphic at $\alp_0$ if $\alp_0$ is a regular point of $\Rep$. Here, $\tau(\alp)$ is defined in \eqref{E:combsection_tau}.

If $\alp_0$ is a regular point of $\Rep$, then by Hartog's theorem, \cite[Th.~2.2.8]{HormanderSCV}, it is enough to show that for every holomorphic curve $\gam:\calO\to U_{\lam,\tet}$, where $\calO$ is a connected open neighborhood of 0 in $\CC$,  the function
\begin{equation}\label{E:f(z)}
    f(z) \ := \ \frac{\rho_{\GG_{[0,\lam]}}(\n_{\gam(z)},g)}{\tau(\gam(z))}
\end{equation}
is complex differentiable at $0$, i.e., there exists $a\in \CC$, such that as $z\to 0$
\[
    f(z)\ =\ f(0)\ +\ z\cdot a\ +\ o(z).
\]

\subsection{Choice of a basis}\label{SS:basis}
We use the notation introduced in \refss{germconnections}. In particular we have a vector bundle $E$, a holomorphic curve $\gamma:\calO\to\Rep$,  and a continuous family of flat connection $\n_{\gam(z)}\ (z\in\calO)$ on $E$ such that for each $z\in \calO$ the monodromy of $\n_{\gam(z)}$ is equal to $\gam(z)$. If $\alp_0$ is a regular point of $\Rep$ then we also assume that
\begin{equation}\label{E:n0+uome2}
    \n_{\gam(z)} \ = \ \n_{\gam(0)} \ + \  z\cdot \ome \ + \ o(z).
\end{equation}
where $o(z)$ is understood in the sense of the Fr\'echet topology.

Let $\Pi_{[0,\lam]}(z)$ ($z\in \calO$) denote the spectral projection of the operator $\B(\n_{\gam(z)},g)^2$, corresponding to the set of eigenvalues of $\B(\n_{\gam(z)},g)^2$, whose absolute value is $\le\lam$. Then it follows from the definition of $U_\lambda$ that  $\Pi_{[0,\lam]}(z)$ depends continuously on $z$. Moreover, in case when $\alp_0$ is a regular point of $\Rep$, $\Pi_{[0,\lam]}(z)$ is complex differentiable in $z$. Hence, in this case there exists a bounded operator
$R$ on $L^2_*(M,E\oplus E)$ such that
\begin{equation}\label{E:Piz-Piz0}
    \Pi_{[0,\lam]}(z) \ = \ \Pi_{[0,\lam]}(0) \ + \ z\, R \ +\ o(z).
\end{equation}

We denote by $\Ome^\b(z)$ the image of $\Pi_{[0,\lam]}(z)$. Recall that we denote the dimension of $M$ by $m=2r-1$. For each $j=0\nek r-1$, fix a basis 
\[
	\bfw_j\ =\ \{w^1_j\nek w^{l_j}_j\}
\]
of \/
$\Ome^j(0)$ and set $\bfw_{m-j}:= \{\GG{}w^1_j\nek \GG{}w^{l_j}_j\}$. To simplify the notation we will write $\bfw_{m-j}= \GG\bfw_j$. Then $\bfw_j$ is a basis for $\Ome^j(0)$ for all $j=0\nek m$.

For each $z\in \calO$, $j=0\nek m$, set
\[
    \bfw_j(z)\ = \ \big\{\,w^1_j(z)\nek w^{l_j}_j(z)\,\big\} \ := \ \big\{\,\Pi_{[0,\lam]}(z)\,w^1_j\nek \Pi_{[0,\lam]}(z)\,w^{l_j}_j\,\big\}.
\]
Since $\Pi_{[0,\lam]}(z)$ depends continuously on $z$, there exists a
neighborhood $\calO'\subset \calO$ of $0$, such that $\bfw_j(z)$ is a basis of $\Ome^j(z)$ for all $z\in \calO'$, $j=0\nek m$. Further, since $\Pi_{[0,\lam]}(z)$ commutes with $\GG$, we obtain
\begin{equation}\label{E:oome=Gamoome}
    \bfw_{m-j}(z)\ =\ \GG\,\bfw_j(z).
\end{equation}
Clearly, $\bfw_j(0)= \bfw_j$ for all $j=0\nek m$.

For each $z\in \calO'$, the space $\Ome^\b(z)$ is a subcomplex of $\big(L^2_*(M,E\oplus E),\DD_{\gam(z)}\big)$. Moreover, the embedding $\Ome^\b(z)\hookrightarrow L^2_*(M,E\oplus E)$  is a quasi-isomorphism. It follows from Theorem~\ref{Thm41} that the cohomology of this compels is canonically isomorphic to 
\[
	H^\b(M,E_{\gam(z)})\oplus H^\b(M,\d M,E_{\gam(z)}).
\]

Let
\[
    \phi_{\Ome^\b(z)}:\, \Det\big(\,\Ome^\b(z)\,\big) \ \longrightarrow \ \Det		\big(\,H^\b(M,E_{\gam(z)})\oplus H^\b(M,\d M,E_{\gam(z)})\,\big)
\]
denote the canonical isomorphism, cf. Section~2.4 of \cite{BraKap:refined-det} . For $z\in \calO'$, let $w(z)\in \Det\big(\,\Ome^\b(z)\,\big)$  be the element determined by the basis
$\bfw_1(z)\nek \bfw_m(z)$ of $\Ome^\b(z)$. More precisely, we introduce
\[
    w_j(z) \ =\  w_j^1(z)\wedge\cdots\wedge{}w_j^{l_j}(z) \ \in \ \Det\big(\,\Ome^j(z)\,\big),
\]
and set
\[
    w(z) \ := \ w_0(z)\otimes{}w_1(z)^{-1}\otimes\cdots\otimes{}w_m(z)^{-1}.
\]
Then, according to Definition~4.3 of \cite{BraKap:refined-det}, it follows from \eqref{E:oome=Gamoome} that, for all $z\in \calO'$, the refined torsion of the complex
$\Ome^\b(z)$ is equal to $\phi_{\Ome^\b(z)}(w(z))$, i.e.,
\begin{equation}\label{E:rho0lam=}
    \rho_{{}_{\GG_{\hskip-1pt{}_{[0,\lam]}}}}(\n_{\gam(z)}) \ = \ \phi_{\Ome^\b(z)}(w(z)).
\end{equation}

\subsection{Reduction to a family of differentials}\label{SS:familyofd}
Using the basis $\bfw_j(z)$ we define the isomorphism
\[
    \psi_j(z):\, \CC^{l_j}  \ {\longrightarrow}\ \Ome^j_{[0,\lam]}(z)
\]
by the formula
\begin{equation}\label{E:psi(z)}
    \psi_j(z)(x_1\nek x_{l_j}) \ := \ \sum_{k=1}^{l_j}\,x_k\,w_j^k(z) \ = \ \sum_{k=1}^{l_j}\,x_k\,\Pi_{[0,\lam]}(z)\,w_j^k.
\end{equation}
We conclude that for each $z\in \calO'$, the complex $\big(\Ome^\b(z),\DD_{\gam(z)}\big)$ is isomorphic to the complex
\begin{equation}\label{E:coordOme0lam}
    \begin{CD}
       (W^\b,d(z)):\quad 0 \ \to \CC^{l_0} @>{d_0(z)}>> \CC^{l_1} @>{d_1(z)}>>\cdots @>{d_{l-1}(z)}>> \CC^{l_m} \ \to \ 0,
    \end{CD}
\end{equation}
where
\begin{equation}\label{E:d(z)}
    d_j(z) \ := \ \psi_{j+1}(z)^{-1}\circ \DD_{\gam(z)}\circ \psi_j(z), \qquad j=0\nek m.
\end{equation}
It follows from \eqref{E:Piz-Piz0} and \eqref{E:psi(z)} that $d_j(z)$ is continuous family of differentials. Moreover, when $\alp_0$ is a regular point of $\Rep$ it is complex differentiable at $0$, i.e., there exists a
$(l_{j+1}\times{}l_j)$-matrix $A$ such that
\[
    d_j(z) \ = \ d_j(0) \ + \ z\,A \ + \ o(z).
\]
Let $\psi(z):= \bigoplus_{j=0}^d\psi_j(z)$. Since $\GG\big(\Ome^j(z)\big)= \Ome^{m-j}(z)$ ($j=0\nek m$), we conclude that $l_j= l_{m-j}$. From
\eqref{E:oome=Gamoome} we obtain that
\begin{equation}\label{E:tilGam}
    \hatGam\ :=\  \psi^{-1}(z)\circ \GG \circ \psi(z)
\end{equation}
is independent of $z\in \calO'$ and
\begin{equation}\label{E:Gamofbasis}
    \hatGam:\, (x_1\nek x_{l_j}) \ \mapsto (x_1\nek x_{l_j}), \qquad j=0\nek m.
\end{equation}
It follows from \eqref{E:d(z)} and \eqref{E:tilGam} that
\begin{equation}\label{E:Gam-tilGam}
    \rho_{{}_{\hatGam}}(z) \ = \ \rho_{{}_{\GG_{\hskip-1pt{}_{[0,\lam]}}}}(\n_{\gam(z)}),
\end{equation}
where $\rho_{{}_{\hatGam}}(z)$ denotes the refined torsion of the finite dimensional complex $(W^\b,d(z))$ corresponding to the chirality operator $\hatGam$.

Let $\phi_{W^\b}(z):\Det(W^\b)\to \Det\big(H^\b(d(z))\big)$ denote the denote the canonical isomorphism of Section~2.4 of \cite{BraKap:refined-det} . The standard bases of $\CC^{l_j}$
($j=0\nek m$) define an element $\tilw\in \Det(W^\b)$. From \eqref{E:Gamofbasis}  and the definition of $\rho_{\GG}(z)$ we
conclude that
\begin{equation}\label{E:rhotilgam}
    \rho_{{}_{\GG}}(z) \ = \ \phi_{W^\b}(z)(\tilw).
\end{equation}

\subsection{The acyclic case}\label{SS:holasyclic}
To illustrate the main idea of the proof let us first consider the case, when both $H^\b(M,E_{\alp_0})$ and 
$H^\b(M,\d M,E_{\alp_0})$ are trivial.  Then there exists a neighborhood $\calO''\subset \calO'$ of $0$ such that $H^\b(M,E_{\gam(z)})= H^\b(M,\d M,E_{\gam(z)})= 0$ for all $z\in \calO''$. Thus the torsion \eqref{E:rhotilgam} is a complex valued function on $\calO''$. To finish the proof of Theorem~\ref{T:holomsection} in this case it remains to show that this function is continuous and, in case when $\alp_0$ is a regular point of $\Rep$,  is complex differentiable at $0$. In view of \eqref{E:d(z)}, this follows from the following
\begin{Lem}\label{L:asyclicanalytic}
Let
\begin{equation*}
 \begin{array}{ll}
    \begin{CD}
       (C^\b,\pa(z)):\quad 0 \ \to \CC^{n\cdot k_0} @>{\pa_0(z)}>> \end{CD}&\begin{CD}\CC^{n\cdot k_1}
       @>{\pa_1(z)}>>\cdots
    \end{CD}\\
    &\begin{CD}
    \cdots@>{\pa_{m-1}(z)}>> \CC^{n\cdot k_m} \ \to \ 0,
    \end{CD}
 \end{array}
\end{equation*}
be a family of acyclic complexes defined for all $z$ in an open set $\calO\subset \CC$. For any $c\in \Det(C^\b)$ the function $z\mapsto \phi_{(C^\b,\pa(z))}(c)$ is continuous if  the differentials $\pa_j(z)$ are continuous, and is complex differentiable at $0$ if $\pa_j(z)$ are complex differentiable at
$0$.
\end{Lem}
\begin{proof}
It is enough to prove the lemma for one particular choice of $c$. To make such  a choice let us fix for each $j=0\nek m$ a complement of
$\IM(\pa_{j-1}(0))$ in $C^j$ and a basis $v_j^1\nek v_j^{l_j}$ of this complement. Since the complex $C^\b$ is acyclic, for all $j=0\nek m$,
the vectors
\begin{equation}\label{E:basisvz0}
    \pa_{j-1}(0)\,v_{j-1}^1\,,\ \ldots\,,\  \pa_{j-1}(0)\,v_{j-1}^{l_{j-1}},\ \,  v_j^1\ ,\,\ldots\,,\  v_j^{l_j}
\end{equation}
form a basis of $C^j$. Let $c\in \Det(C^\b)$ be the element defined by these bases. Then, for all $z$ close enough to $0$  and for all
$j=0\nek m$,
\begin{equation}\label{E:basisvz}
    \pa_{j-1}(z)\,v_{j-1}^1\,,\ \ldots\,,\  \pa_{j-1}(z)\,v_{j-1}^{l_{j-1}},\ \, v_j^1\,,\ \ldots\,,\   v_j^{l_j}
\end{equation}
is also a basis of $C^j$. Let $A_j(z)\ (j=0\nek m)$ denote the non-degenerate matrix transforming the basis \eqref{E:basisvz} to the basis
\eqref{E:basisvz0}. Then, by the definition of the isomorphism $\phi_{(C^\b,\pa(z))}$, cf. \S2.4 of \cite{BraKap:refined-det} ,
\begin{equation}\label{E:phi(z)(c)}
    \phi_{(C^\b,\pa(z))}(c)\ = \ (-1)^{\calN(C^\b)}\ \prod_{j=0}^m\,\Det(A(z))^{(-1)^j},
\end{equation}
where $\calN(C^\b)$ is the integer defined in formula (2.15) of \cite{BraKap:refined-det}  which is independent of $z$. Clearly, the matrix
valued functions $A_j(z)$ and, hence, their determinants are continuous if  the differentials $\pa_j(z)$ are continuous, and are complex differentiable at $0$ if $\pa_j(z)$ are complex differentiable at
$0$. Thus, so is the function $z\mapsto
\phi_{(C^\b,\pa(z))}(c)$.
\end{proof}

\subsection{Sketch of the proof of Theorem~\ref{T:holomsection} in the non-acyclic case}\label{SS:sketchnonacyclic}
We now turn to the proof of  Theorem~\ref{T:holomsection} in the general case. In this subsection we sketch the main ideas of the proof. It is enough to show that the function
\begin{equation}\label{E:f(z)2}
    f(z) \ := \ \frac{\rho_{{}_{\GG_{\hskip-1pt{}_{[0,\lam]}}}}(\n_{\gam(z)},g)}{\tau(\gam(z))}
\end{equation}
continuous and, if $\alp_0$ is a regular point of $\Rep$, is complex differentiable at $0$. Here $\tau$ is the map \eqref{E:combsection_tau}.
To see this we consider the de Rham integration maps
\begin{equation}
\begin{split}
    &J_z^{\max}:\, \Ome^\b_{\max}(M,E_{\gam(z)}) \ \longrightarrow \ C^\b(K,\gam(z)), \\
    &J_z^{\min}:\, \Ome^\b_{\min}(M,E_{\gam(z)}) \ \longrightarrow \ C^\b(K, K',\gam(z)).
\end{split}
\end{equation}
where the cochain complexes $C^\b(K,\gam(z))$ and $C^\b(K,K',\gam(z))$ are defined in \S \ref{S:detlinebundle}.

The de Rham integration map of $E$-valued differential forms is defined using a trivialization of $E$ over each cell $e_j$, and, hence, it depends on the flat
connection $\n_{\gam(z)}$, cf. below. More precisely, in the neighborhood of any cell $e$ of $K$, a differential form $f\in \Ome^\b(M,E_{\gam(z)})$
can be written in the form $ f= \sum_{j=1}^n f_j \otimes v_j$,
where $f_j \in \Ome^\b(M)$ is a complex-valued differential form and $v_j$ is a $\nabla_{\gamma(z)}$-flat section of 
$E_{\gam(z)}$ for $j=1,..,n$. The de Rham integration map is then defined by
\[
J^{\max}_z f(e) := \sum_{j=1}^n \left(\int_e  f_j \right) v_j.
\]
If $f\in \Omega^\b_{\min}(M,E_{\gamma(z)})$ then $J^{\max}_z f$ is a well-defined element of 
the relative cochain complex $C^\b(K,K',\gam(z))$ and we denote the corresponding map by $J^{\min}_z$. 
Both maps descend to isomorphisms on cohomology \cite[\S 4]{RS}. We write 
\begin{equation}
  \begin{split}
	&\Omega^\b_d(M,E_{\alp}) := \Omega^\b_{\min}(M,E_{\alp}) \oplus \Omega^\b_{\max}(M,E_{\alp}), \\
	&C^\b_d(K, \alp) := C^\b(K, K', \alp) \oplus C^\b(K,\alp).
  \end{split}
\end{equation} 
Hence we obtain a quasi-isomorphism 
\begin{align}
J_z := J^{\max}_z \oplus J^{\min}_z: \Omega^\b_d(M,E_{\gamma(z)}) \to C^\b_d(K, \gamma(z)).
\end{align}

The trivialization $T_e(z):E|_e \to \C^n \times e$ induced by the connection $\nabla_{\gamma(z)}$ is continuous. Moreover, if $\alp_0$ is a regular point of $\Rep$, it is complex differentiable
at $z=0$ by \eqref{E:n(z)}. Hence so is the de Rham integration map $J_z$. 

We then consider the restriction $J_z|_{\Ome^\b(z)}$ of $J_z$ to the finite dimensional complex $\Ome^\b(z)$ and study the cone complex $\Cone^\b(J_z|_{\Ome^\b(z)})$ of the map $J_z$. This is a finite dimensional acyclic complex
with a fixed basis, obtained from the bases of $\Ome^\b(z)$, defined in \S \ref{SS:basis}, 
and $C^\b_d(K,\gam(z))$. The torsion of this complex is equal to $f(z)$. An
application of \refl{asyclicanalytic} to this complex proves Theorem~\ref{T:holomsection}.

In the definition of the integration map $J_z$ we have to take into account the fact that the vector bundles $E_{\gam(z)}$ and $E=E_{\gam(0)}$ are
isomorphic but not equal. The integration map $J_z$, cf. \refss{combin}, is a map from $\Ome^\b_d(M,E_{\gam(z)})$ to the cochain complex $C^\b_d(K,\gam(z))$,
which is not equal to the complex $C^\b_d(K,\gam(0))$. Fix an Euler structure on $M$. It defines an  isomorphism between the complexes $C^\b_d(K,\gam(0))$ and $C^\b_d(K,\gam(z))$ which depends on $z$. The study of this isomorphism, which is conducted in \refss{tilC-C}, is important for the understanding of the properties of $J_z$. In particular, it is used to show that in a certain sense $J_z$ is complex differentiable at $0$, which implies that the cone complex $\Cone^\b(J_z)$ satisfies the conditions of \refl{asyclicanalytic}.

\subsection{The cochain complex of the bundle $E$}\label{SS:combin}
Fix a CW-decompo\-si\-tion $K= \{e_1\nek e_N\}$ of $M$. For each $j=1\nek N$ choose a point $x_j\in e_j$ and let $E_{x_j}$ denote the fiber of $E$
over $x_j$. Then the cochain complex $(C^\b_d(K, \gam(z)), \partial_\b)$ may be naturally identified with the complex 
$(C^\b_d(K, E), \partial'_\b(z))$ where the $z$-dependence is now fully encoded in the differentials.
We use the prime in the notation of the differentials $\pa^\prime_j$ in order to distinguish them from the differentials of the cochain complex
$C^\b_d(K,\gam(z))$. The differentials $\pa^\prime_j(z)$ are continuous. Moreover, if $\alp_0$ is a regular point if $\Rep$ then it follows from \eqref{E:n0+uome2} that $\pa^\prime_j(z)$ are complex differentiable at $0$, i.e., there exist linear maps
$ a_j:\,C^j_d(K,E) \longrightarrow\ C^{j+1}_d(K,E)$, s.t.
\[
    \pa^\prime_j(z) \ = \ \pa^\prime_j(0)\ + \ z\cdot a_j\ + \ o(z), \qquad j=1\nek m-1.
\]

\subsection{Relationship with the complex $C^\b_d(K,\gam(z))$}\label{SS:tilC-C}
Recall that for each $z\in \calO'$ the monodromy representation of $\n_{\gam(z)}$ is equal to $\gam(z)$. Let $\pi:\tilM\to M$ denote the universal cover of $M$ and let $\tilE= \pi^*E$ denote the pull-back of the bundle $E$ to $\tilM$. Recall that in \refss{germconnections} we fixed a point $x_*\in M$. Let $\tilx_*\in \tilM$ be a lift of $x_*$ to $\tilM$ and fix a basis of the fiber $\tilE_{\tilx_*}$ of $\tilE$ over $x_*$. Then, for each $z\in \calO'$, the flat connection $\n_{\gam(z)}$ identifies $\tilE$ with the product $\tilM\times \CC^n$. 

Recall from   \refss{Turaevtorsion} that the choice of the  Euler structure $\eps$ also fixes the lifts $\tile_j \ (j=1\nek N)$  of the cell
$e_j$ fixed in \refss{Turaevtorsion}. Let $\tilx_j\in \tile_j$ be the lift of $x_j\in e_j$. Then the trivialization of $\tilE$ defines a continuous in $z$ family of isomorphisms
\[
    S_{z,j}:\,E_{x_j} \ \simeq \ \tilE_{\tilx_j} \ \to \ \CC^n, \qquad j=1\nek N, \ z\in \calO'.
\]
The isomorphisms $S_{z,j}$ depend on the trivialization of $\tilE$, i.e., on the connection $\n_{\gam(z)}$. The direct sum $S_z= \bigoplus_jS_{z,j}$ defines
an isomorphism $S_z:C^\b_d(K,E)\to C^\b_d(K,\gam(z))$. If $\alp_0$ is a regular point of $\Rep$, it follows from \eqref{E:n(z)} that
$S_z$ is complex differentiable at $0$, i.e. for some linear map $s$
\[
    S_z \ = \ S_{0}\ + \ z\cdot s\ + \ o(z).
\]

Finally, we consider the morphism of complexes
\begin{equation}\label{E:tilintegr}
    I_z \ := \ S^{-1}_z\circ J_z\circ \psi(z):\, W^\b\ \to C^\b:=C^\b_d(K,E), \ z\in \calO'.
\end{equation}
This map is complex differentiable at $z=0$ and induces an isomorphism of cohomology.

\subsection{The cone complex}\label{SS:cone}
The cone complex $\Cone^\b(I_z)$ of the map $I_z$ is given by the sequence of vector spaces
\[
    \Cone^j(I_z) \ := \ W^j\oplus C^{j-1}_d\big(K,E) \ \simeq \ \CC^{l_j}\oplus \CC^{n\cdot k_{j-1}}, 
\]
$j=0\nek m$, with differentials
\[
    \hatpa_j(z) \ = \ \left(\,\begin{array}{cc} d_j(z)&0\\I_{z,j}&\pa'(\gam(z))\end{array}\,\right),
\]
where $I_{z,j}$ denotes the restriction of $I_z$ to $W^j$. This is a family of acyclic complexes with differentials $\hatpa_j(z)$, which are continuous. If $\alp_0$ is a regular point of $\Rep$ then $\hatpa_j(z)$ are also
complex differentiable at $0$. The standard bases of $\CC^{l_j}\oplus \CC^{n\cdot k_{j-1}}$ define an element $c\in \Det(\Cone^\b(I_z))$ which
is independent of $z\in \calO'$. Using the canonical isomorphism of Section~2.4 of \cite{BraKap:refined-det}, we hence obtain for each $z\in \calO'$ the number $\phi_{\Cone^\b(I_z)}(c)\in \CC\backslash\{0\}$. From the discussion in \refss{sketchnonacyclic} it follows that this number is equal to the ratio \eqref{E:f(z)}. Hence, to finish the proof of the Theorem~\ref{T:holomsection} it remains to show that the function $z\mapsto \phi_{\Cone^\b(I_z)}(c)$ is continuous and is complex differentiable at $0$ if $\alp_0$ is a regular point of $\Rep$. This follows immediately from \refl{asyclicanalytic}.

\section{Gluing formula for refined analytic torsion}

Let $(M,g)$ be a closed oriented Riemannian manifold and $(N,g^N)$ a separating 
hypersurface, such that $M=M_1 \cup_N M_2$. The metric $g$ restricts to Riemannian metrics on 
the two compact components $M_1$ and $M_2$. Assume that $g$ is product in an open tubular neighborhood of $N$.

\subsection{The transmission complex $\Ome^\b(M_1\texttt{\#} M_2,E)$}\label{SS:complexM1M2}
A given representation $\alp \in \textup{Rep}(\pi_1(M), \C^n)$ induces a connection $\D_\alp$ on a vector bundle $E\equiv E_\alp$,
which restricts to well-defined connections on $M_{1,2}$. 
We denote by  $\rho_i:= \rho(\D_\alp, M_i)$ the refined analytic torsions on $M_i, i=1,2$ and by  $\rho=\rho(\D_\alpha, M)$ for the refined analytic torsion on $M$, cf. \eqref{rat-definition}. Let $\iota_j: N \hookrightarrow M_j$ denote the obvious inclusions, $j=1,2$.

We define the transmission subcomplex of $\Omega^\b(M_1, E_\alp) \oplus \Omega^\b(M_2, E_\alp)$ by 
specifying \emph{transmission} boundary conditions
\begin{equation*}
\begin{split}
\Omega^\b(M_1 \texttt{\#} M_2, E) &:= \{(\w_1, \w_2) \in \Omega^\b(M_1, E_\alp) \oplus \Omega^\b(M_2, E_\alp) \mid \iota^*_1\w_1 = \iota^*_2\w_2\}, \\
\nabla_\alp(\w_1, \w_2) &:= (\nabla_\alp\w_1, \nabla_\alp\w_2).
\end{split}
\end{equation*}
This defines a complex with eigenforms of the corresponding Laplacian given by the eigenforms of the 
Hodge-Laplacian on $(\Omega^\b(M,E_\alp), \nabla_\alp)$, cf. \cite[Theorem 5.2]{Vertman10}. 
In particular their de Rham cohomologies coincide. 

\subsection{The fusion map}\label{SS:gluingfusion}
The splitting
$M=M_1\cup_N M_2$ now gives rise to short exact sequences of the associated complexes
\begin{equation}\label{SES}
\begin{split}
&0 \to \Omega^\b_{\min}(M_1, E) \xrightarrow{\alp_1} \Omega^\b(M_1 \texttt{\#} M_2, E) \xrightarrow{\beta_2} \Omega^\b_{\max}(M_2, E) \to 0, \\
&0 \to \Omega^\b_{\min}(M_2, E) \xrightarrow{\alp_2} \Omega^\b(M_1 \texttt{\#} M_2, E) \xrightarrow{\beta_1} \Omega^\b_{\max}(M_1, E) \to 0,
\end{split}
\end{equation}
where $\alp_1(\w)=(\w, 0), \alp_2(\w)=(0,\w)$ and $\beta_j(\w_1, \w_2) =\w_j, j=1,2$.
The corresponding long exact sequences in cohomology yield canonical isomorphisms, cf. \cite{Vertman10}
\begin{equation*}
\begin{split}
& \Phi_1(\alp): \Det (H^\b(M_1, N,E)) \otimes \Det (H^\b(M_2, E)) \to \Det (H^\b(M,E)),\\
& \Phi_2(\alp): \Det (H^\b(M_2, N,E)) \otimes \Det (H^\b(M_1, E)) \to \Det (H^\b(M,E)).
\end{split}
\end{equation*}
The fusion isomorphisms, cf. \cite[(2.18)]{BraKap:refined-det} provide canonical identifications 
\begin{equation*}
\begin{split}
& \mu_1: \Det (H^\b(M_1, N,E)) \otimes \Det (H^\b(M_1, E)) \to \Det (H^\b(\domm_1, \DD_1)),\\
& \mu_2: \Det (H^\b(M_2, N,E)) \otimes \Det (H^\b(M_2, E)) \to \Det (H^\b(\domm_2, \DD_2)),\\
& \mu: \Det (H^\b(M, E)) \otimes \Det (H^\b(M, E)) \to \Det (H^\b(\domm, \DD)),
\end{split}
\end{equation*}
where the Hilbert complexes $(\domm_j, \DD_j)$ are defined in Definition \ref{domain}, 
with the lower index $j$ referring to the underlying manifold $M_j, j=1,2$. The Hilbert 
complex $(\domm, \DD)$ is defined over $M$. We put
\begin{multline*}
\Phi \equiv \Phi(\alp) \: = \mu \circ (\Phi_1(\alp) \otimes \Phi_2(\alp)) \circ (\mu_1^{-1} \otimes \mu_2^{-1})
:\\ \Det (H^\b(\domm_1, \DD_1)) \otimes \Det (H^\b(\domm_2, \DD_2))\ \to \Det (H^\b(\domm, \DD)).
\end{multline*}
If $\alp \in \Rep$ is unitary, \cite[Theorem 10.6]{Vertman10} asserts
\begin{align}\label{app-4}
\Phi \left( \rho_1 \otimes \rho_2\right) = K \cdot \rho,
\end{align}
where  
\begin{equation}\label{E:K=}
	K \equiv K(\alp)\\ =\ \sigma \,2^{\chi(N)} e^{ \pi i (
		\eta(\B_\alpha, M) - \eta(\B_\alpha, M_1) - \eta(\B_\alpha, M_2))}. 
\end{equation}
Here $\chi(N)$ stands for the Euler characteristic of  $(N, E_\alp\restriction N)$ and 
the $\eta$-invariants $\eta(\B_\alpha, M)$, $\eta(\B_\alpha, M_1)$ and $\eta(\B_\alpha, M_2)$ are defined 
in terms of the even parts of the corresponding odd signature operators. The sign $\sigma$ is determined by formula (8.4) of \cite{Vertman10}. Note that the sign depends  on the dimensions of various cohomology spaces and is related to the sign convention used in defining the fusion isomorphism of determinant lines, see \cite[\S2]{BraKap:refined-det}  for a detail discussion of the sign conventions. Note also that  the Euler characteristic  $\chi(N)$ depends only on the rank $n$ of the representation $\alp$ but not on the particular choice of $\alp\in \Rep$. 

The $\eta$-invariants $\eta(\B_\alpha, M)$, $\eta(\B_\alpha, M_1)$ and $\eta(\B_\alpha, M_2)$ are not necessarily continuous functions of $\alp$: they have integer jumps when some eigenvalues of the odd signature operator cross zero. Hence, $K^2(\alp)$ is a continuous function of $\alp$. Moreover, below we show that $K^2$ extends to a weakly holomorphic function on the space $\Rep$ of representations. 

\subsection{The gluing formula for some non-unitary representations}\label{SS:extendnonunitary}
The main result of this section is the following extension of \eqref{app-4} to some class of non-unitary representations.
\begin{Thm}\label{refined-gluing-theorem}
Let $(M,g)$ be a closed oriented Riemannian manifold of odd dimension, and $N$ a separating 
hypersurface such that $M=M_1 \cup_N M_2$, and $g$ is product in an open tubular neighborhood of $N$. Assume  that  $\calC\subset \Rep$ is a connected component and $\alp_0\subset \calC$ is a unitary representation which is a regular point of the complex analytic set $\calC$. For any representation  $\alp\in\Rep$ denote by $\rho(\alp)$ and $\rho_j(\alp)$ the refined analytic 
torsions on $M$ and $M_j, j=1,2$. Then for any $\alp\in \calC$ we have
\begin{align*}
	\Phi \left( \rho_1(\alp) \otimes \rho_2(\alp)\right)\ =\ \pm\,K(\alp) \cdot \rho(\alp).
\end{align*}
\end{Thm}
The rest of this section is occupied with the proof of Theorem~\ref{refined-gluing-theorem}, which is based on an analytic continuation technique, cf. \cite[\S6.4]{BraVer}. First we need the following proposition.

\begin{Prop}\label{K}
$K(\alp)^2$ is a weakly holomorphic  function on the complex analytic space $\Rep$.
\end{Prop}
\begin{proof}
We need to show that $\exp\big(2i\pi\eta(\B_\alpha, M_i))$, $(i=1,2)$ and \linebreak$\exp\big(2i\pi\eta(\B_\alpha, M)\big)$ are weekly holomorphic functions of $\alpha$.

Denote by $\B_{\alp,j}$ ($j=1,2$) the odd signature operator of the complex $(\domm_j,\DD_j)$ and by $\B_\alp$ the odd signature operator of the comples $(\domm,\DD)$. With this notation we have 
\[
	\eta(\B_\alpha, M_j)\ = \ \eta(\B_{\alpha,j}),\qquad
	\eta(\B_\alpha, M)\ = \ \eta(\B_{\alpha}).
\]

Fix $\alp_0\in \Rep$ and  a number $\lam>0$ such that the spectra of  the operators $\B_{\alp_0}$ and $\B_{\alp_0,j}$ ($j=1,2$) do not intersect the circle $\{z\in\CC:\, |z|= \lambda\}$. There exists a neighborhood $U_\lam$ of $\alp_0$ in $\Rep$ such that for all $\alp\in U_\lam$ the spectra of the operators $\B_{\alp}$ and $\B_{\alp,j}$ ($j=1,2$) do not intersect the circle $\{z\in\CC:\, |z|= \lambda\}$. Then (cf. formula (4-1) of \cite{Vertman09})
\begin{equation}\label{E:B-Blam}
	\eta(\B_{\alp,j})\ - \ \eta(\B_{\alpha,j}^{(\lam,\infty)})\
	 \in \frac12\,\ZZ,
	\qquad
	\eta(\B_\alpha)\ - \ \eta(\B_{\alpha}^{(\lam,\infty)})\
	 \in \ \frac12\,\ZZ.
\end{equation}
Notice that the functions  $\eta(\B_{\alpha,j}^{(\lam,\infty)})$ and $\eta(\B_{\alpha}^{(\lam,\infty)})$ are continuous on $U_\lam$. The functions $\eta(\B_{\alp,j})$ and $\eta(\B_\alpha)$ are not necessarily continuous, but might have integer jumps. Hence, it follows from \eqref{E:B-Blam} that 
$\eta(\B_{\alp,j})-\eta(\B_{\alpha,j}^{(\lam,\infty)})$ and $\eta(\B_\alpha)- \eta(\B_{\alpha}^{(\lam,\infty)})$ are constants modulo $\ZZ$. We conclude that  
\[
	\exp\big(2i\pi\eta(\B_{\alp,j})\ - \ 
	\exp\big(2i\pi\eta(\B_{\alp,j}^{(\lam,\infty)}),
	\qquad j=1,2,
\]
and  
\[
	\exp\big(2i\pi\eta(\B_{\alp})\ - \ 
	\exp\big(2i\pi\eta(\B_{\alp}^{(\lam,\infty)})
\]
are constant functions on $U_\lam$. Hence, it suffices to show that the functions $\exp\big(2i\pi\eta(\B_{\alp,j}^{(\lam,\infty)})$, $(i=1,2)$ and $\exp\big(2i\pi\eta(\B_\alp^{(\lam,\infty)})\big)$ are weakly holomorphic in a neighborhood of $\alp_0$. 

Let $\theta$ be an Agmon angle for the operators $B_{\alp_0}$ and $B_{\alp_0,j}$. Then there exists a neighborhood $U_{\lam,\theta}\subset U_\lam$ of $\alp_0$, such that for all $\alp\in U_{\lam,\theta}$  $\theta$ is an Agmon angle for $B_{\alp}$ and $B_{\alp,j}$. By \cite[(4.6)]{Vertman09} we obtain
\[
	\Det'_{gr,\theta}(\B_{\alp,j, \textup{even}}^{(\lambda,\infty)})\
	 = \
	 \exp \Big(\,
	\xi_\lam(\alp,M_j) - i\pi \xi'_\lam(\alp,M_j) -
	i \pi \eta(\B_\alpha^{(\lam,\infty)}, M_j)\, \Big),
\]
where we have set 
\begin{align*}
	&\xi_\lam(\alp,M_j)\ :=\  
	\frac{1}{2} \sum_{k=0}^m (-1)^k \cdot k \cdot \left. \frac{d}{ds} 
	\right|_{s=0} 
	\zeta_{2\theta}\big(s, \B^2_{\alp,j} \restriction 
	   \domm_{j,(\lam,\infty)}^k\big), \\
	&\xi'_\lam(\alp,M_j)\ :=\
	\frac{1}{2} \sum_{k=0}^m (-1)^k \cdot k \cdot 
	\zeta_{2\theta}\big(s=0, \B^2_{\alp,j} \restriction 
		\domm_{j,(\lam,\infty)}^k\big).
\end{align*}

The graded determinant $\Det'_{gr,\theta}(\B_{\alp,j, \textup{even}})$ is weakly holomorphic in $U_{\lam,\theta}$ by Corollary~\ref{C:analytic}.
The fact that  $\exp\big(2\xi(\alp,M_j)\big)$ and \linebreak$\exp\big(2i\pi\xi'(\alp,M_j)\big)$ are weakly holomorphic follows similarly from the variational formula
\begin{multline}\notag
	\left.\frac{d}{dt} \right|_{t=0}\!\!\!
	\zeta_{2\theta}\big(s, \B^2_{\kappa(t)} \restriction 
	\domm_{j,(\lam,\infty)}^k\,\big)
	\\ =\ 
    -s \cdot \textup{Tr} \, 
    \Big(\,\frac{d}{dt}\Big|_{t=0} \,\B^2_{\kappa(t)}
    \restriction 
	\domm_{j,(\lam,\infty)}^k\,\Big)  		   
    \big(\,\B^2_{\kappa(0)}\restriction 
	\domm_{j,(\lam,\infty)}^k\,\big)^{-s-1}.
\end{multline}
As a consequence, $\exp\big(2i\pi\eta(\B_{\alpha,j}^{(\lam,\infty)}))$\ ($j=1,2$) are weakly holomorphic. Similarly, on proves that $\exp\big(2i\pi\eta(\B_\alpha^{(\lam,\infty)})\big)$ is weakly holomorphic.
\end{proof}

\subsection{An analytic continuation}\label{SS:continuation}
The set of unitary representations is the fixed point set of the anti-holomorphic involution 
\[
	\tau:\,\Rep\ \to \Rep, \qquad \tau:\,\alp\ \mapsto \ \alp',
\]
where $\alp'$ denotes the representation dual to $\alp$. Hence, it is a totally  real submanifold of $\Rep$ whose real dimension is equal to $\dim_\CC\calC$, cf. \cite[Proposition~3]{Ho04}. In particular there is a holomorphic coordinates system $(z_1\nek z_r)$ near $\alp_0$ such that the unitary representations form a {\em real neighborhood} of $\alp_0$, i.e. the set $\IM{z_1}=\ldots =\IM{z_r}=0$. Therefore, cf. \cite[p.~21]{ShabatSeveralVariables}, we obtain the following 

\begin{Prop}\label{P:continuation}
 If two holomorphic functions coincide on the set of unitary representations they also coincide on $\calC$. 
\end{Prop}

\subsection{The proof of Theorem~\ref{refined-gluing-theorem}}\label{SS:prrefined-gluing-theorem}

By Theorem \ref{T:holomsection} the refined analytic torsions $\rho_1,\rho_2$ and $\rho$ 
define holomorphic sections on the corresponding determinant line bundles.  The canonical isomorphism $\Phi$ defines a bilinear map between holomorphic determinant line bundles and hence maps holomorphic sections to holomorphic sections. Let us denote by $f(\alp)$ the unique complex valued function of $\alp$ such that 
\[
	\Phi \left( \rho_1 \otimes \rho_2\right) = f(\alp) \cdot \rho,
\]
holds for all $\alp\in \calC$. Since both $\Phi \left( \rho_1 \otimes \rho_2\right)$ and $\rho$ are holomorphic sections of the determinant line bundle, $f(\alp)$ is a holomorphic function. 

It follows from \eqref{app-4} that 
\begin{equation}\label{E:f2=K2}
	K(\alp)^2\ = \ f(\alp)^2
\end{equation}
for all unitary representations in $\calC$. By Proposition~\ref{K}, $K(\alp)^2$ is a holomorphic function. Hence, we obtain from  Proposition~\ref{P:continuation} that \eqref{E:f2=K2} holds for all $\alp\in \calC$. 
\hfill$\square$

\begin{Rem}\label{R:ref-usual}
The gluing formula for refined analytic torsion may be used to prove a gluing result for the Ray-Singer 
torsion norm on connected components of the representation variety that contain a unitary point.
\end{Rem}

\section{Gluing formula for Ray-Singer analytic torsion}

We continue in the previously outlined setup of a closed oriented Riemannian 
manifold $(M,g)$ and a separating hypersurface $(N,g^N)$, such that $M=M_1\cup_N M_2$.
Consider a representation $\alp\in \Rep$ and the corresponding flat vector bundle $(E,\nabla)$.
Fix a Hermitian metric $h_0$ on $E$, of product type near $N$. Even if the metric structures $(g,h_0)$ are product near $N$, 
the connection $\nabla$ need not be product near $N$, so that the resulting Laplacian does not 
have a product structure near the separating hypersurface and hence a gluing theorem for 
Ray-Singer analytic torsion cannot be obtained from the results of \cite{Lue}, \cite{Vis} and \cite{Les}. 

Before we proceed, let us make some chronological remarks on that topic.
Vishik \cite{Vis} was first to prove the gluing formula for analytic torsion given a unitary representation 
without using the theorem of Cheeger \cite{Che} and M\"uller \cite{Mue1}. Though it was not explicitly stated in 
\cite{Vis}, the assumption of a unitary representation is obsolete once a connection is in temporal gauge\footnote{ The notion
and properties of temporal gauge are recalled in the Appendix.} near $N$. Under the assumption of 
temporal gauge and product metric structures, the Hodge Laplacian is of product type near $N$ and the Vishik's 
argument goes through. Recently, Lesch \cite{Les} provided an excellent discussion of the gluing formula 
for possibly non-compact spaces, extending the result of Vishik, and stating clearly that the proof requires only product metric structures
and a connection in temporal gauge near the cut, rather than unitariness of the representation. 

A quite general proof of the gluing formula for general representations and without assuming product
metric structures, was provided by Br\"uning and Ma \cite{BruMa:gluing}. They derive a gluing formula by 
relating the Ray-Singer and the Milnor torsions, in odd and also in the even-dimensional case. In this section 
we present a different proof of their result in odd dimensions, as a consequence of \cite{Les} and the 
Br\"uning-Ma anomaly formula in \cite{BruMa}. To make the main idea of our alternative argument clear, 
we restrict ourselves to the situation, when $\partial M =\varnothing$.

\begin{Prop}\label{BM-hidden}
Consider two Hermitian metrics $h_0,h_1$ on a fixed flat vector bundle $(E,\nabla)$
over a compact oriented odd-dimensional Riemannian manifold $(K,g), j=1,2$. If $\partial K \neq \varnothing$, 
assume that $h_0, h_1$ coincide over $Y=\partial K$. Fix either relative or absolute boundary 
conditions at $Y$ for the Hodge Laplacian and denote by $\| \cdot \|^{\textup{RS}}_{(g,h_j)}, j=0,1,$ 
the corresponding Ray-Singer analytic torsion norms. Then \footnote{\ There is no product structure assumption on $g$ and $h_j, j=0,1$.}
\[
\| \cdot \|^{\textup{RS}}_{(g,h_0)} = \| \cdot \|^{\textup{RS}}_{(g,h_1)}.
\]
\end{Prop}
\begin{proof}
The metric anomaly, identified in \cite{BruMa} is expressed in terms of 
the Levi-Civita connections $\nabla^{TK}$ and $\nabla^{TY}$ on $K$ and its boundary $Y$, 
the respective representatives $e(TK,\nabla^{TK}), e(TY, \nabla{TY})$ of the 
Euler classes of $TK, TY$ in Chern-Weil theory, and the quotient 
$\|\cdot \|_{h_0} / \|\cdot \|_{h_1}$ between the metrics on $\det E$, 
induced by $h_0$ and $h_1$.

Since $\dim K$ is odd, $e(TK,\nabla^{TK})=0$. If $h_0, h_1$ coincide over $Y$, 
$\log \|\cdot \|_{h_0} / \|\cdot \|_{h_1}=0$ over $Y$, so that the statement 
follows from \cite[Theorem 0.1, (0.5)]{BruMa}.
\end{proof}

The main idea now is the reduction to the setup of a connection in temporal gauge near $N$. 
A connection $\nabla$ is in temporal gauge in an open neighborhood $\mathscr{U}=(-\epsilon, \epsilon) \times N$ 
of $N$, if $\nabla = \pi^*\nabla^N$ for some flat connection $\nabla^N$ on $E_N$, where
$\pi:\mathscr{U} \to N$ is the natural projection onto the second factor. Proposition \ref{t-g} below asserts
that in fact every connection is gauge equivalent to a connection in temporal gauge.

We denote the corresponding gauge transformation by $\gamma$. 
The gauge transformed connection is given by $\nabla_\gamma = \gamma \nabla \gamma^{-1}$.
We set for any $u, v\in \Gamma(M,E)$
\[
h_\gamma (u,v) := h_0 (\gamma u , \gamma v).
\]
This defines a new Hermitian metric on $E$ that coincides with $h$ over $N$, since 
$\gamma$ acts as identity over $N$. By construction, $\gamma$ induces an isometry
\[
\gamma: L^2_*(M,E;g,h_\gamma) \to L^2_*(M,E;g,h_0).
\]
The following theorem is a result by Lesch \cite{Les}, cf. Vishik \cite{Vis}.
\begin{Thm}\label{lesch}
Let $(M,g)$ be a closed oriented Riemannian manifold of odd dimension, and $N$ a separating 
hypersurface such that $M=M_1 \cup_N M_2$, and $g$ is product in an open tubular neighborhood of $N$.
Consider flat Hermitian vector bundle $(E,\nabla_\gamma, h_0)$. Denote by $\| \cdot \|^{\textup{RS}}$ and 
$\| \cdot \|^{\textup{RS}}_{M_i}$ the corresponding Ray-Singer norms on $M$ and $M_i,i=1,2$, 
respectively. The Ray-Singer norm on $M_1$ is defined with respect to relative boundary conditions, while 
on $M_2$ we pose absolute boundary conditions. 
Then
\begin{align*}
\log \frac{\| \cdot \|^{\textup{RS}}_M}{\Phi_1 \left(\| \cdot\|^{\textup{RS}}_{M_1} 
\otimes \| \cdot \|^{\textup{RS}}_{M_2}\right)} 
= \frac{1}{2} \chi(N) \log 2.
\end{align*}
\end{Thm}

\begin{Cor}
Let the Ray-Singer torsion norms $\| \cdot \|^{\textup{RS}}$ and 
$\| \cdot \|^{\textup{RS}}_{M_i}$ be defined with respect to the flat connection $\nabla$
and any Hermitian metric $h$ on $E$. Then 
\begin{align*}
\log \frac{\| \cdot \|^{\textup{RS}}_M}{\Phi_1 \left(\| \cdot\|^{\textup{RS}}_{M_1} 
\otimes \| \cdot \|^{\textup{RS}}_{M_2}\right)} 
= \frac{1}{2} \chi(N) \log 2.
\end{align*}
\end{Cor}

\begin{proof} 
The Laplacian $\Delta_\gamma$ on $E_\gamma = (E,\nabla_\gamma, h_0)$ is related to the Laplacian 
$\Delta$ on $E_\alp =(E,\nabla, h_\gamma)$ by the unitary transformation $\gamma$
with $\Delta_\gamma = \gamma \circ \Delta \circ \gamma^{-1}$. Hence, $\gamma$
induces a map between the harmonic forms of $\Delta$ and $\Delta_\gamma$, and
hence also between the corresponding determinant lines, which we also denote by $\gamma$.
By construction, we find 
\begin{equation}
\| \gamma (\cdot) \|^{\textup{RS}}_{E_\gamma} = \| \cdot \|^{\textup{RS}}_{E_\alp},
\end{equation}
where we indicate the dependence on the vector bundle by the subindex 
and omit the reference to the underlying manifold, since the relation holds both on $M$ and $M_i, i=1,2$.
The isometric identification $\gamma$ commutes with the maps in the
\eqref{SES}. Hence $\Phi_1 \circ (\gamma \otimes \gamma) = \gamma \circ \Phi_1$, 
and by Theorem \ref{lesch}, we find
\begin{align}\label{BM-gluing-formula}
\log \frac{\| \cdot \|^{\textup{RS}}_{(M,E_\alp)}}{\Phi_1 \left(\| \cdot\|^{\textup{RS}}_{(M_1,E_\alp)} 
\otimes \| \cdot \|^{\textup{RS}}_{(M_2,E_\alp)}\right)} 
= \frac{1}{2} \chi(N) \log 2.
\end{align}
This is a gluing theorem for any complex representation $\alp$, possibly non-unitary.
A priori, however, this relation holds for the specific Hermitian metric $h_\gamma$ on $(E,\nabla_\alp)$.
Since $\gamma \restriction N = \textup{id}$, the metrics
$h_0, h_\gamma$ coincide over $N$ and hence, by Proposition \ref{BM-hidden} 
the gluing formula \eqref{BM-gluing-formula} holds for any representation $\alp$
and any choice of a Hermitian metric $h$ on $(E,\nabla_\alp)$. 
\end{proof}

The Riemannian metric $g$ is still assumed to be product near $N$, so that variation of $g$ 
leads to additional anomaly terms, cf.  \cite[Theorem 0.1]{BruMa}.

\section{Appendix: Temporal Gauge Transformation}\label{appendix-gauge}

In this section we recall the notion of a connection in temporal gauge and review 
some properties of these connections, cf. \cite{Vertman10}. In particular, we show that any 
flat connection is gauge equivalent to a connection in temporal gauge, cf. Proposition \ref{t-g}.
Consider a closed oriented Riemannian manifold $(M,g^M)$ of dimension $m$ and a vector bundle 
$E$ with structure group $G\subset GL(n,\C)$. Denote the principal $G$-bundle associated to $E$ by $P$,
where $G$ acts on $P$ from the right.

Consider a hypersurface $N\subset M$ and its collar neighborhood $U\cong (-\epsilon,\epsilon)\times N$.
We view the restrictions $P|_U,P|_{N}$ as $G$-bundles, where the structure group can possibly be reduced to a subgroup of $G$. 
Let $\pi: (-\epsilon, \epsilon)\times \partial X \to \partial X$ be the natural projection onto the second component. 
Then $E|_U\cong \pi^*E|_{N}$ and for the associated principal bundles $P|_U\cong \pi^*P|_{N}\xrightarrow{f} P|_{N}$,
where $f$ is the principal bundle homomorphism, covering $\pi$, with the associated homomorphism of the structure groups 
being the identity automorphism. 

\begin{Def}
We call a flat connection $\w$ on $P$ a connection in temporal gauge over $U$, 
if there exists a flat connection $\w_N$ on $P_N$ such that $\w|_U=f^*\w_{N}$ 
over the collar neighborhood $U$. Similar condition defines a covariant derivative in 
temporal gauge.
\end{Def}

We now explain the notion of temporal gauge in local terms.
Let $\w_{N}$ denote a flat connection one-form on $P|_{N}$. Then $\w_U:=f^*\w_{N}$ 
gives a connection one-form on $P|_U$ which is flat again. 
Let $\{\widetilde{U}_{\beta},\widetilde{\Phi}_{\beta}\}_\beta$ be a set of local trivializations for $P|_{N}$. 
Then $P|_U\cong \pi^*P|_{N}$ is trivialized over the local neighborhoods 
$U_{\beta}:=(-\epsilon, \epsilon)\times \widetilde{U}_{\beta}$ with the induced trivializations $\Phi_{\beta}$. 
For any $y \in \widetilde{U}_{\beta}$, normal variable $x \in (-\epsilon, \epsilon)$ and for $e\in G$ being the 
identity matrix we put
\begin{align*}
\widetilde{s}_{\beta}(y):=\widetilde{\Phi}_{\beta}^{-1}(y,e), \quad s_{\beta}(x,y):=\Phi_{\beta}^{-1}((x,y),e).
\end{align*}
These local sections define local representations for $\w_U$ and $\w_{N}$
\begin{align*}
\widetilde{\w}^{\beta}:=\widetilde{s}^*_{\beta}\w_{N}\in \Omega^1(\widetilde{U}_{\beta},\mathcal{G}), \\
\w^{\beta}:=s^*_{\beta}\w_U\in \Omega^1(U_{\beta},\mathcal{G}),
\end{align*}
where $\mathcal{G}$ denotes the Lie algebra of $G$. Consider local coordinates $y=(y_1,..,y_{m-1})$ on $\widetilde{U}_{\beta}$. Then 
\begin{equation}\label{temp-stern}
\begin{split}
&\widetilde{\w}^{\beta}=\sum\limits_{i=1}^n \w^{\beta}_i(y)dy_i, \quad 
\w^{\beta}=\w^{\beta}_0(x,y)dx+\sum\limits_{i=1}^n \w^{\beta}_i(x,y)dy_i, \\
&\textup{with} \quad \w^{\beta}_0\equiv 0, \ \textup{and} \ \w^{\beta}_i(x,y)\equiv \w^{\beta}_i(y). 
\end{split}
\end{equation}

\begin{Prop}\label{t-g}
Any flat connection on the principal bundle $P$ is gauge equivalent to a flat connection in temporal gauge.
\end{Prop} 
\begin{proof}
By a partition of unity argument it suffices to discuss the problem locally over $U_{\beta}$. 
Let $\w$ be a flat connection on $P|_U$. Let $\gamma\in \textup{Aut}(P|_U)$ be any gauge transformation on $P|_U$
and $\gamma^\beta$ the corresponding local representation. Denote the gauge transform of $\w$ under $\gamma$ by $\w_\gamma$. 
The local $\mathcal{G}$-valued one-forms $\w^{\beta}, \w_\gamma^\beta$ are related in correspondence to the transformation law 
of connections by
$$
\w^{\beta}_\gamma=(\gamma^\beta)^{-1}\cdot \w^{\beta}\cdot \gamma^\beta+(\gamma^\beta)^{-1}d\gamma^\beta,
$$
where the action $\cdot$ is the multiplication of matrices ($G\subset GL(n,\C)$), after evaluation at a local vector field and 
a base point in $U_{\beta}$. 

The local one form $\w^{\beta}$ can be written as
$$
\w^{\beta}=\w^{\beta}_0(x,y)dx+\sum\limits_{i=1}^n w^{\beta}_i(x,y)dy_i.
$$
Our task is to identify the correct gauge transformation $\gamma$, so that $\w$ is temporal gauge, 
cf. \eqref{temp-stern}. For this reason we consider the following initial 
value problem with parameter $y\in \widetilde{U}_{\beta}$
\begin{equation} \label{ODE} 
\begin{split}
\partial_x\gamma^\beta(x,y)&=-\w^{\beta}_0(x,y)\gamma^\beta(x,y), \\
\gamma^\beta(0,y)&=\one \in GL(n,\C).
\end{split}
\end{equation}
The solution to \eqref{ODE} is given by an integral curve of the time dependent vector field 
$V^{\beta}_{x,y}$ on $G$, parametrized by $x\in (-\epsilon, \epsilon)$, 
such that for any $u\in G$ 
$$
V^{\beta}_{x,y}u:=-(R_u)_* \w^{\beta}_0(x,y)= -\w^{\beta}_0(x,y)\cdot u,
$$
where $R_u$ is the right multiplication on $\gamma$ and the second equality follows from the fact that 
$G\subset GL(n,\C)$ is a matrix Lie group. The corresponding integral curve $\gamma^{\beta}(x, y)$ 
with $\gamma^{\beta}(0, y)=\one\in G$ is $G$-valued and the unique solution to \eqref{ODE}.

We now compute for the gauge transformed connection $\w_\gamma$
\begin{align*}
\w^{\beta}_\gamma&=(\gamma^\beta)^{-1}\cdot \w^{\beta}\cdot \gamma^\beta+(\gamma^\beta)^{-1}d\gamma^\beta\\
&=(\gamma^\beta)^{-1}\cdot \w^{\beta}_0\cdot \gamma^\beta dx+\sum\limits_{i=1}^n (\gamma^\beta)^{-1}\cdot 
\w^{\beta}_i\cdot \gamma^\beta dy_i \\ &+ (\gamma^\beta)^{-1}\partial_x \gamma^\beta dx +\sum\limits_{i=1}^n 
(\gamma^\beta)^{-1}\partial_{y_i} \gamma^\beta dy_i \\
&=\sum\limits_{i=1}^n (\gamma^\beta)^{-1}\cdot \w^{\beta}_i\cdot \gamma^\beta dy_i + 
\sum\limits_{i=1}^n (\gamma^\beta)^{-1}\partial_{y_i} \gamma^\beta dy_i.
\end{align*} 
where in the last equality we cancelled two summands due to $\gamma^\beta$ being the solution to \eqref{ODE}. 
We now use the fact that $\w$ is a flat connection. A gauge transformation preserves flatness, so $\w_\gamma$ is flat again. Put
$$
\w^{\beta}_\gamma=\w^{\beta}_{\gamma,0}(x,y)dx+\sum\limits_{i=1}^n \w^{\beta}_{\gamma,i}(x,y)dy_i,
$$ 
where by the previous calculation 
$$
\w^{\beta}_{\gamma,0}\equiv 0, \quad \w^{\beta}_{\gamma,i}\equiv (\gamma^\beta)^{-1}
\cdot \w^{\beta}_i\cdot \gamma^\beta+ (\gamma^\beta)^{-1}\partial_{y_i} \gamma^\beta.
$$
Flatness of $\w_\gamma$ implies 
$$
\partial_x\w^{\beta}_{\gamma,i}(x,y)=\partial_{y_i}\w^{\beta}_{\gamma,0}(x,y)=0.
$$
Hence the gauge transformed connection is indeed in temporal gauge. This completes the proof.
\end{proof}

\def\cprime{$'$} \def\cprime{$'$} \newcommand{\noop}[1]{} \def\cprime{$'$}
\providecommand{\bysame}{\leavevmode\hbox to3em{\hrulefill}\thinspace}
\providecommand{\MR}{\relax\ifhmode\unskip\space\fi MR }
\providecommand{\MRhref}[2]{%
  \href{http://www.ams.org/mathscinet-getitem?mr=#1}{#2}
}
\providecommand{\href}[2]{#2}

\end{document}